\documentclass[12pt]{article}

\usepackage{epic}
\usepackage{eepic}
\usepackage{epsf} 
\usepackage{ amssymb, amscd}
\usepackage{amsmath}
\usepackage{color}
\numberwithin{equation}{section}
\usepackage{epsfig}
\usepackage{graphicx}
\usepackage{tikz}
\usetikzlibrary{matrix}

\def\cA            {{\mathcal{A}}}
\def\cB            {{\mathcal{B}}}
\def\cC            {{\mathcal{C}}}
\def\cD            {{\mathcal{D}}}

\def\cF            {{\mathcal{F}}}
\def\cG            {{\mathcal{G}}}
\def\cH            {{\mathcal{H}}}

\def\cM            {{\mathcal{M}}}

\def\cV            {{\mathcal{V}}}

\def\cZ            {{\mathcal{Z}}}

\def\bbC           {\mathbb{C}}

\def\bbM           {\mathbb{M}}

\def\bbR           {\mathbb{R}}
\def\bbT           {\mathbb{T}}
\def\bbZ           {\mathbb{Z}}

\def\MXN           {{}_M {\mathcal{X}}_N}
\def\MXM           {{}_M {\mathcal{X}}_M}

\def\MXMp          {{}_M^{} {\mathcal{X}}_M^+}
\def\MXMm          {{}_M^{} {\mathcal{X}}_M^-}
\def\MXMpm         {{}_M^{} {\mathcal{X}}_M^\pm}

\def\NXN           {{}_N {\mathcal{X}}_N}

\def\NXM           {{}_N {\mathcal{X}}_M}

\def\sdprod{{\times\!\vrule height5pt depth0pt width0.4pt\,}}

\title{Reconstruction and Local Extensions for Twisted Group Doubles, and Permutation Orbifolds}
\author{
{\sc David E.\ Evans}\\
 {\footnotesize School of Mathematics, Cardiff University,}\\
 {\footnotesize Senghennydd Road, Cardiff CF24 4AG, Wales, U.K.}\\
 {\footnotesize e-mail: {\tt EvansDE@cf.ac.uk}}\\ \\
 {\sc Terry  Gannon }\\
 {\footnotesize Department of Mathematics, University of Alberta,}\\
{\footnotesize Edmonton, Alberta, Canada T6G 2G1}\\
{\footnotesize e-mail: {\tt tgannon@math.ualberta.ca}} }

\begin{document}
\maketitle

\begin{abstract} 
We prove the first nontrivial reconstruction theorem for modular tensor categories:  the category associated to any twisted Drinfeld
double of any finite group, can be realised as the representation category of a completely rational conformal net. We also show that any twisted double of a solvable  group is the category of modules of  a completely rational  vertex operator algebra. In the process of doing this, we identify the 3-cocycle twist for permutation orbifolds of holomorphic conformal nets: unexpectedly, it can be nontrivial, and depends on the value of the central charge modulo 24. In  addition, we determine the branching coefficients of all possible {local (conformal)} extensions of any finite group orbifold of
holomorphic conformal nets, and identify their modular tensor categories. All statements also apply to vertex operator algebras, provided the conjecture holds that finite group orbifolds of holomorphic VOAs are rational, with a category of modules given by a twisted group double.  

\end{abstract}

{\footnotesize
\tableofcontents
}

\section{Introduction}

The finite-dimensional complex representations of a finite group $G$ form a semi-simple rigid tensor category Rep($G$). It has finitely many simples (one of which is the tensor unit), and the Hom-spaces are finite-dimensional vector spaces. Such categories are called  \textit{fusion}. Moreover, the tensor product is symmetric: the obvious isomorphism $c_{\rho,\phi}:\rho\otimes\phi\rightarrow
\phi\otimes\rho$ satisfies $c_{\rho,\phi}\circ c_{\phi,\rho}=id_{\phi\otimes\rho}$.  A version of Tannaka--Krein duality due to Deligne (see e.g. Corollary 9.9.25 of \cite{EGNO}) says that any symmetric fusion category is Rep($G$) (though possibly with the braidings $c_{\rho,\phi}$ twisted by an order 2 element $z$ in the centre of $G$) --- this is called \textit{reconstruction} --- and $G,z$ are unique up to equivalence.

In section 2.1 we define complete rationality for vertex operator algebras (VOAs) and
conformal nets of factors.
The modules of completely rational VOAs $\cV$ and conformal nets  $\cA$ also form fusion categories Mod$(\cV)$ and Rep$(\cA)$ respectively, but now the braiding isomorphisms $c_{M,N}$ are  as far as possible from being symmetric. Such categories are called \textit{modular tensor categories (MTC)}. One can ask whether every MTC is equivalent to some Mod$(\cV)$ or Rep($\cA$) (this would be reconstruction), and if so, to what extent this realisation is unique. Little is known about either question, but it has been
conjectured that any \textit{unitary} MTC can be realised as both Mod$(\cV)$ and Rep$(\cA)$. If so, the realisation will be far from unique, and it is far from clear what plays the role of the central involution $z$ here.

An MTC is a remarkable object. It carries representations of the modular group SL$_2(\bbZ)$ and indeed every surface mapping class group, gives link invariants for every closed 3-manifold, etc. The trivial MTC is the category Vec$_\bbC$ of vector spaces, with only one simple object. 
Reconstruction is easy for Vec$_\bbC$: the VOAs and conformal nets with only one simple module are called \textit{holomorphic}. There are infinitely many holomorphic VOAs and conformal nets (e.g.
the Monstrous Moonshine VOA), and classifying all of them is hopeless. 

The easiest classes of unitary MTC, as well as completely rational VOAs and conformal nets, are associated to even positive-definite lattices $L$. The simples are in natural bijection
with the cosets $L^*/L$ of the dual lattice by $L$, so self-dual $L$ are holomorphic. All simples in these MTC 
are invertible (such categories are called \textit{pointed}); it is not difficult to
show that any unitary pointed MTC is the category of modules of a lattice VOA and lattice conformal net
\cite{EG7}.

The class of unitary MTC we consider in this paper are the \textit{twisted Drinfeld doubles $\cD_\omega(G)$ of finite groups,} where $[\omega]\in H^3(G;\bbT)$. Because of their importance and ease of construction, they have a long history, starting with their introduction by Dijkgraaf--Witten \cite{DW}. The simple objects are pairs $[g,\chi]$ where $g\in G$  and $\chi$ is a projective irrep of the centraliser $C_G(g)$, whose multiplier 2-cocycle comes from $\omega$. These group doubles have a natural interpretation in terms of equivariant $G\times G$ bundles, and this geometric picture plays an important role in our arguments.    The importance of twisted group doubles lies in the orbifold construction, one of the fundamental constructions in the theory.  If $G$ is a finite group of automorphisms of a  VOA or conformal net $\cA$, then we define the orbifold $\cV^G$ or $\cA^G$ to be  the space of points fixed by  $G$.  When $\cA$ is holomorphic, the category Rep$(\cA^G)$ will be $\cD_\omega(G)$ for some twist $\omega$ \cite{Xu},\cite{Kir} (the analogue is conjectured to hold for VOAs).

The easiest and best-studied class of orbifolds are
the \textit{permutation orbifolds}.
Let $\cA$ (or $\cV$) be as above, and let $G$ be any subgroup of some symmetric group $S_k$. The permutation orbifold of $\cA$ by $G$ is $(\cA\otimes\cdots\otimes\cA)^G$ where $G$ acts by permuting the $k$ copies of $\cA$. When $\cA$ is holomorphic, so will be $\cA\otimes\cdots\otimes\cA$, and  the  MTC of the permutation orbifold will be  $\cD_\omega(G)$ for some twist $[\omega]\in H^3(G;\bbT)$. It has long been believed (see e.g. Conjecture 6.3 in \cite{Mug1}, and section IV(B) in \cite{Dav2}) that the class $[\omega]$ must be trivial, because this construction doesn't see any structure in $\cA$. In Theorem 2 we show that $[\omega]$ is \textit{not} always trivial, though it depends only on the central charge $c$ (a basic numerical invariant) of the net. This shows that the MTC Rep$((\cA^{\otimes k})^G)$ of permutation orbifolds  does not depend only on $G$ and the MTC Rep$(\cA)$ (which here is trivial). In Theorem 4(a) we prove the same result  for VOAs when $G$ is solvable and $c$ is a multiple of 24. The possible nontriviality of $\omega$ can be seen quite easily using characters, as we do at the end of section 3.2. Our proof of Theorem 2 relies on work by Nakaoka \cite{Nak} on the cohomology of the symmetric group.

We use this to prove reconstruction
 (Theorem 3)  for the twisted group doubles: for any group $G$ and twist $[\omega]\in H^3(G;\bbT)$, $\cD_\omega(G)$ is Rep$(\cA^G)$ for some holomorphic net $\cA$. This is the first nontrivial class of examples of reconstruction.  Our proof also works for VOAs, provided Conjecture 1 below holds. In particular, in Theorem 4(b) we establish reconstruction by VOAs for any twisted double of a solvable group (no conjectures are assumed). An important part of our reconstruction proofs uses the 8-term restriction-inflation exact sequence in cohomology.

A big part of group theory is to understand, and exploit, the subgroups $H$ in a group $G$. The relation between Mod($G$) and Mod($H$) is described by the induction and restriction functors. The analogous notion for VOAs and conformal nets is \textit{local} or \textit{conformal extension}: just as the category
Mod$(H)$ is generally smaller the smaller $H$ gets, the category Mod$(\cV)$ gets smaller the larger
$\cV$ gets. In  categorical language, extensions  $\cV^e$ of $\cV$ correspond to \textit{module categories}
of Mod$(\cV)$ of local extension type, also called \textit{type 1}. Again, induction and restriction functors describe things. For example, the conformal extensions of MTC, VOAs and conformal nets of a lattice $L$, correspond to 
even lattices $L'$ containing $L$ but of the same dimension, so $L\subseteq L'\subseteq L^{\prime *}\subseteq L^*$. For example, restriction takes the $\cA(L')$-irrep corresponding to the coset
$v+L'\in L^{\prime*}/L'$ to $v+L'\subset L^*/L$, i.e. to a sum of $|L'/L|$ irreps for $\cA(L)$. 

To prove reconstruction we need explicit control of the conformal extensions $(\cA^G)^e$ of the holomorphic orbifolds $\cA^G$ (and of $\cV^G$). Equivalently, we need control on all type 1 module categories of any twisted double $\cD_\omega(G)$. We do that in Theorem 1. We explicitly give the restriction and
induction functors (for example the branching coefficients of $(\cA^G)^e$-modules to $\cA$-modules), as well as the MTC of the extension $(\cA^G)^e$ (it itself is always  a twisted  group double).  These extensions are built from  two basic classes: the obvious class corresponds to the orbifold $\cA^K$ by some subgroup $K\le G$, and the MTC of the extension is simply $\cD_\omega(K)$; the other  is more difficult to describe, and has MTC $\cD_{\omega'}(G/N)$ for some normal subgroup $N$ of $G$. Every conformal extension is an extension (possibly trivial) from the first class, followed by one (possibly trivial) from the second.  Determining the twist $\omega'$ in the second class is quite subtle. Our proof builds on Davydov--Simmons \cite{Dav3} and  the thesis \cite{Jon2} of Vaughan Jones. (Incidentally, we describe the remaining module categories of $\cD_\omega(G)$, e.g. those of type 2, in \cite{EG8}.) 

For example, consider Rep$(\cA^G)\cong\cD_1(G)$ for any $G$. The extreme example of the first class
is $K=1$: the resulting extension recovers the holomorphic net $\cA$ itself, with branching coefficients
Res$(\cA)=\sum_\chi\mathrm{dim}\,\chi\,[1,\chi]$.  The extreme example of the second class is $\overline{G}=1$, also
a holomorphic net, with branching coefficients $\sum_g[g,1]$.

The twisted doubles are absolutely no less fundamental than the untwisted ones --- e.g.  Monstrous Moonshine concerns the orbifold of the Moonshine VOA by its full automorphism group, the Monster finite simple group $\bbM$, and the corresponding $\cD_\omega(\bbM)$ has order-24 twist $[\omega]$ \cite{JF}. The significance of the 3-cocycle $\omega$ is obstruction: the more nontrivial  $\omega$ is, the fewer conformal extensions $\cA^G$  and
$\cV^G$ have, and the fewer (and larger) are their irreps. This is illustrated in Table 1 below.

This paper falls into the sequence of  papers \cite{EG1},\cite{EG3},\cite{Ev},\cite{EG7},\cite{EG8} on the $K$-theory of loop groups and finite groups. Though it can be read independently, we were led to many of our arguments by this picture,  made explicit for finite groups in \cite{Ev}, and we use bundles throughout.

There are two obvious extensions of this work, and we intend to address both in the near future. One is reconstruction for the doubles of other fusion categories associated to finite groups. For example, the weakly group theoretical fusion categories conjecturally exhaust all fusion categories with integer global dimension. We would expect these to be constructed as nested orbifolds starting from a holomorphic net (or VOA). The general case is out of reach at present, but we have completed reconstruction for the  Tambara--Yamagami categories (reconstruction of the easier half, that of even rank, has also been announced in \cite{Bis}). The importance of this is that these can be regarded as the building blocks of doubles of fusion categories such as the Haagerup and other quadratic categories.

The second related extension of this work is to arbitrary finite group orbifolds of conformal nets (or VOAs) whose MTC
of representations is pointed (such as the lattice theories). The resulting MTC are also weakly group theoretical. In recent work, Mason and Ng \cite{MN} have conjectured the explicit form of these categories, as well as the associated quasi-Hopf algebra. The Verlinde rings have a natural $K$-theoretic interpretation, so seem to be a natural generalisation of the finite group doubles. The challenge is reconstruction for the MTC of \cite{MN},  as say orbifolds of say lattice theories.

Section 2 reviews the twisted Drinfeld double $\cD_\omega(G)$ of finite groups, as well as notions of bundles  and module categories for finite groups. Section 3 states our results and gives several examples, and section 4 supplies the proofs.

\section{Finite group doubles}

\subsection{Module categories, alpha induction, and all that}

For background on {modular tensor category} (MTC) and fusion categories, see e.g.\ \cite{EGNO}.
 Let $\Phi$ denote the (finite) set of {isomorphism classes $\lambda$ of} irreducible objects in $\cC$; we call these \textit{sectors}.
The Grothendieck ring of $\cC$ is called the \textit{Verlinde} or \textit{fusion ring} Ver, {and has 
basis $\Phi$}.  

We occasionally allude to the \textit{modular data} of an MTC, which is a unitary representation of
the modular group SL$_2(\bbZ)$ on the complexification $\bbC\otimes_\bbZ\mathrm{Ver}$.
{Since SL$_2(\bbZ)$ is generated by  $\left({0\atop 1}{ -1\atop  \,0}\right)$ and $\left({1\ 1\atop 0\ 1}
\right)$, this representation is uniquely determined by the matrices $S,T\in M_{\Phi\times\Phi}(\bbC)$ corresponding to those generators.} $T$ is a diagonal matrix, {whereas $S$ determines the tensor product structure constants of Ver through Verlinde's formula.} 

For the basic theory of vertex operator algebras (VOAs), see e.g. \cite{LL}. By a \textit{completely rational VOA}, we mean $C_2$-cofinite, regular, simple, self-dual and of CFT-type. Then the category Mod$(\cV)$ of modules of a completely rational VOA is a MTC \cite{Hua}.  For the basic theory of local conformal nets, see e.g. \cite{Kaw}. By a \textit{completely rational conformal net}, we mean one with the split property and finite $\mu$-index. Then the category Rep$(\cA)$ of representations of a completely rational conformal net is a unitary MTC \cite{KLM}. It is expected that the theory of sufficiently nice (e.g. completely rational and fully unitary) VOAs and their modules should be naturally equivalent to that of sufficiently nice (e.g. completely rational) conformal nets --- see \cite{CKLW} for recent progress. 

This paper focuses on  the twisted group doubles $\cD_\omega(G)$, where $G$ is a finite group and $\omega\in Z^3(G;\bbT)$. These unitary  MTC are described next subsection. They arise in the theory of VOAs and conformal nets as follows.  A completely rational VOA $\cV$ or conformal net $\cA$ is called \textit{holomorphic} if its representation theory  is trivial, i.e. Mod$(\cV)$ resp.\ Rep$(\cA)$ is Vec$_\bbC$. If $G$ is a finite group of automorphisms of a VOA $\cV$ or conformal net $\cA$ (so by definition  $ G$ acts faithfully on the underlying spaces $\cV$ resp.\ $\cH$), by the  \textit{orbifold}  $\cV^G$ resp.\ $\cA^G$ we mean the associated VOA resp.\ conformal net structure on the space of fixed points of all $g\in G$ --- see e.g.\ \cite{LL},\cite{Xu} for details. If $\cA$ is completely rational then so is  $\cA^G$ \cite{Xu}; the analogue is conjectured for VOAs (see e.g. \cite{CM} for recent progress). When $\cA$ is holomorphic, the category Rep$(\cA^G)$ will be $\cD_\omega(G)$ for some twist $\omega$ \cite{Kir} (the analogue holds for $\cV$ when $\cV^G$ is completely rational).

Any completely rational VOA or conformal net comes with a representation of the Virasoro algebra. In this reprsentation, the (normalised) Virasoro central term $C$ is sent to a multiple of the identity. The \textit{central charge}  $c$ of a conformal net or VOA is that numerical factor.
The corresponding MTC determines $c$ only up to mod 8. In particular, the central charge of a 
holomorphic net or VOA is a multiple of 8. The central charge enters Theorem 2 below.

A \textit{conformal} or \textit{local extension} $\cV^e$ resp.\ $\cA^e$ of a completely rational VOA $\cV$ or conformal net $\cA$ is a completely rational VOA or conformal net containing 
$\cV$ resp.\ $\cA$ but with the same central charge, so that  any $\cV^e$-module or $\cA^e$-representation restricts to  Mod$(\cV)$ resp.\ Rep$(\cA)$. These restrictions are called \textit{branching rules}. Such extensions correspond to commutative symmetric special Frobenius  algebras $A$ in the MTC $\cC=\mathrm{Mod}(\cV)$ or Rep$(\cA)$ as follows (for VOAs this is developed in \cite{KO,HKL,CKM}, and the same arguments work for conformal nets). Here, $A$ is the restriction of $\cV^e$ or $\cA^e$ to
$\cV$ resp.\ $\cA$; to be an algebra it must have a multiplication $\mu\in\mathrm{End}_{\cC}(A)$; to be commutative this must satisfy $\mu\circ c_{A,A}=\mu$ where  $c_{A,A}$ is the braiding. For the definition of symmetric special Frobenius (which we don't need), see e.g. \cite{FFSS}.  A left $A$-module is a pair
$(V,\mu_V)$, where $V\in\cC$ and $\mu_V:A\otimes V\rightarrow V$ is the $\cC$-morphism corresponding to multiplication. These form a 
 fusion category  Mod$_\cC(A)$, where the tensor product $V\otimes_A W$ is some summand of $V\otimes_{\cC} W$. Alpha-induction is the functor $\alpha:\cC\rightarrow \mathrm{Mod}_\cC(A)$ defined on objects by $\alpha(V)=A\otimes V$ and
$\mu_{\alpha(V)}=\mu\otimes id$ where $\mu$ is the multiplication in $A$, and sigma-restriction is the forgetful functor
Res$:\mathrm{Rep}\,A\rightarrow \cC$ sending $(V,\mu_V)$ to $V$. {Alpha-induction is a tensor functor, and sigma-restriction is its adjoint (i.e. Frobenius reciprocity holds between them)}.  The category
Mod$_\cC(A)$ is not in general braided; the MTC Mod$(\cV^e)$ or Rep$(\cA^e)$ is identified with the full subcategory Mod$_\cC{}^0(A)$ 
of dyslectic {or local} objects. The simple
objects in Mod$_\cC{}^0(A)$ consist of the simple $(V,\mu_V)\in\mathrm{Mod}_\cC(A)$ with twist $\theta_V\in\bbC\,id$.

More generally, we have the notion of \textit{module category} \cite{Ost1}. 
The MTC $\cC$ only describes (part of) the \textit{chiral} data of the rational CFT;  its module category  captures its boundary data, defect lines, spaces of conformal
blocks in arbitrary genus, etc (see e.g. \cite{FFSS}).  A module category $\cM$ over an MTC $\cC$
is a bifunctor $\otimes:\cC\times\cM\rightarrow\cM$ together with compatible associativity and unit
isomorphisms. These satisfy the usual pentagon and triangle identities, corresponding to different ways to 
identify $((X\otimes Y)\otimes Z)\otimes M\cong X\otimes(Y\otimes(Z\otimes M))$ and $(X\otimes 1)\otimes M\cong
X\otimes M$ for all objects $X,Y,Z$ of $\cC$ and $M\in\cM$. There are obvious notions of equivalence and direct sums
of module categories, and of indecomposable module categories. The Main Theorem of \cite{Ost1} says that each indecomposable module category of  $\cC$ is the category Mod$_\cC(A)$ of right modules in $\cC$ of a  symmetric special Frobenius algebra $A\in\cC$. When $A$ is not commutative, Mod$_\cC(A)$ will no longer be a tensor category. For this reason one also considers the  $A$-$A$-bimodules in $\cC$.
These form a tensor category, called the \textit{full system} $\cC^*_{\cM}=Fun_{\cC}(\cM,\cM)$, a (typically nonbraided) fusion category Morita-equivalent  to (i.e. with the same double as) $\cC^{op}$. There are two ways, called
$\alpha_\pm(V)$,  to make the left-module $A\otimes V$ into a bimodule: we get the right-module structure through either a braiding $c_{A,V}$ or inverse braiding $c_{V,A}^{-1}$. Then both alpha-inductions $\alpha_\pm$ are tensor functors.

This categorical picture was abstracted from the older subfactor picture \cite{LR,O,xu,BE1,BE4,BEK1,BEK2}. The objects in the various categories are unital $*$-homomorphisms
$\rho,\rho':A\rightarrow B$ between type III factors $A,B$. The tensor product is composition; equivalences (called \textit{sectors}) $[\lambda]$ and sums can be defined. Let  $\NXN$ be a system of endomorphisms on
 a factor $N$ which realises the Grothendieck ring  Ver of $\cC$; we write $\Sigma(\NXN)$ for the set of formal sums of those endomorphisms.  
Let $\iota:N\rightarrow M$ be the inclusion $N\subset M$ of factors, and let $\overline{\iota}:M\rightarrow N$ be its conjugate.
 We require   $\NXN$ to be braided, and  {the \textit{dual canonical endomorphism}} 
 $\theta=\overline{\iota}\iota$ to be expressible by $\NXN$.  Using the braiding  or its opposite, 
 we can lift an endomorphism $\lambda\in\NXN$
of $N$ to one of $M$ in two ways, namely the {alpha-inductions}  $\alpha_{\pm }(\lambda)$.  
The induced systems $\MXMpm=\alpha_\pm(\NXN)$
 generate  the {full system} $\MXM$.    Sigma-restriction is $\overline{\iota}\beta\iota$.
By  $\MXN$ we 
mean all irreducibles  appearing in any  $[\lambda\overline{\iota}]$ for
$\lambda\in\NXN$. The \textit{nimrep}, precursor to the notion of module category, is the $\NXN$ action on  {$\Sigma(\NXM)$}, given by left composition. The algebra $A$ 
is $\theta$, with the $Q$-system structure providing the Frobenius algebra structure.

 Given sectors $\lambda,\mu\in\NXN$ define $\langle \lambda,\mu\rangle=\delta_{\lambda,\mu}$, and extend bilinearly to $\Sigma(\NXN)$;
 then  $\langle \lambda,\mu\rangle=\mathrm{dim \,Hom}_\cC(\lambda,\mu)$.  
 The matrix defined by\begin{equation}\label{modinvalpha}\cZ_{\lambda,\mu}:=\langle \alpha_+(\lambda),\alpha_-(\mu)\rangle\end{equation}
is the \textit{modular invariant} associated to the system, since it commutes with $S$ and $T$.
The {modular invariant} helps describe how the full CFT is built from chiral data.

When the Frobenius algebra $A$ is commutative, Mod$_\cC(A)$ is a module category over $\cC$,
with additive functor from the Deligne product $\cC\boxtimes \mathrm{Mod}_\cC(A)$ to $\mathrm{Mod}_\cC(A)$ sending $(V,X)$ to $\alpha(V)\otimes_A X$. Such module categories are called \textit{type 1} (or extension type); as we know they correspond to local extensions. In this paper we focus on the type 1 module categories of twisted group doubles. The  modular invariant of {type 1} module categories is block diagonal. A module category is called \textit{type 2} (or automorphism type)  if its modular invariant is a permutation matrix. 
Every module category is a combination of two type 1's and a type 2, in the sense we explain next.

Given a full CFT, it is natural to speak of the maximally extended left and right chiral algebras, called the \textit{type 1 parents}.  In the subfactor picture \cite{BE4}, these correspond to intermediate subfactors
$N\subset M_{\pm}\subset M$ so that $N\subset M_\pm$ are type 1, with dual canonical endomorphisms
$[\overline{\iota}_+\iota_+]=\sum_{\lambda\in\Phi}\cZ_{\lambda,1}[\lambda]$ and  
$[\overline{\iota}_-\iota_-]=\sum_{\lambda\in\Phi}\cZ_{1,\lambda}[\lambda]$ respectively, coming
from the first row and column of $\cZ$. {In the category language \cite{FFRS}, the type 1 parents are called the left and right centre of the corresponding module category or algebra.} The type 1 parents have equivalent MTC; they can both be canonically identified with a subsystem of the full system $\MXM$ generated by the intersection {$\MXMp\cap\MXMm$},  called the \textit{neutral system} $\MXM^0$.
Write $b_\pm$ for the branching rules, written in matrix form,  from the type 1 parents  ${}_{M_\pm}\mathcal{X}_{M_\pm}^0$ to {$\cC=\NXN$}.
Then \eqref{modinvalpha} becomes in matrix form
\begin{equation}\cZ=b_+\sigma b_-^t\,,\label{modinvsigma}\end{equation}
 where  $\sigma$ is a permutation matrix corresponding to the composition of the canonical identifications mentioned above. When the module category  is type 1, the intermediate factors satisfy $M_+=M_-=M$ and $\sigma=id$ so \eqref{modinvsigma} collapses to $\cZ=b_+b_+^t$.
 For this reason, the branching rules are fundamental and  the modular invariant will henceforth be ignored in this paper.

\subsection{The category $\cD_\omega(G)$}

Throughout this paper,  let Irr$_c(G)$ denote the set
of all {isomorphism classes of projective} irreps of $G$ {with 2-cocycle 
multiplier $c$}, and $R_G$  the character ring of $G$. 
We write $g^h$ for $h^{-1}gh$, ${}^hg$ for $hgh^{-1}$ and cl$_G(g)$ for the conjugacy class
$\{g^h:h\in G\}$. For any subgroup $K\le G$, write $C_K(g)=\{k\in K:g^k=g\}$; e.g. $C_G(g)$ is the centraliser. We write $\Delta_G$ for the diagonal subgroup $\{(g,g)\}$ of $G^2=G\times G$.

Projective representations of finite groups arise naturally in this theory. We assume the reader is familiar
with their basic theory. A standard reference is \cite{Karp}; the paper  \cite{Che} treats the theory
in parallel to that of linear (a.k.a. true) representations, which is also the philosophy we adopt.

Recall the 2- and 3-cocycle conditions for finite group cohomology:
\begin{align}\label{2coc}\psi(x,y)\,\psi(xy,z)=&\,\psi(y,z)\,\psi(x,yz)\,,\\
\label{3coc}\omega(g,h,k)\,\omega(g,hk,l)\,\omega(h,k,l)=&\,\omega(gh,k,l)\,\omega(g,h,kl)\,.\end{align}
For us these always  take values in the unit circle $\bbT\subset\bbC$. We call $\psi\in Z^2(G;\bbT)$  \textit{normalised} when all $\psi(g,1)=\psi(1,h)=1$, and  $\omega\in Z^3(G;\bbT)$  normalised when
all $\omega(1,h,k)=\omega(g,1,k)=\omega(g,h,1)=1$. Given any functions
$f:G\rightarrow \bbT$ and $F:G^2\rightarrow\bbT$, we get a 2-coboundary by ${\delta f}(g,h)=f(gh){f(g)^*}{f(h)^*}$ 
{(where as always we denote complex conjugation with `$*$')} and a 3-coboundary by ${\delta F}(g,h,k)=\frac{F(g,hk)}{F(g,h)}\frac{F(h,k)}{F(gh,k)}$; then {$H^i(G;\bbT)$ is the group of $i$-cocycles quotient that of the $i$-coboundaries.} 
Given a 3-cocycle $\omega\in Z^3(G;\bbT)$, define  $c^\omega_g=c_g$ by
\begin{equation} \label{2cocycle}c_g(h_1,h_2)={\omega(g,h_1,h_2)}\,{\omega(h_1,h_2,g^{h_1h_2})}\,\omega(h_1,g^{h_1},h_2)^*\end{equation}
(a special case of the slant product $H_m\times H^n\rightarrow H^{n-m}$).
Then {for all $g,h_1,h_2,h_3\in G$,} {$\omega$ normalised implies $c_1(h_1,h_2)=c_g(1,h_2)=c_g(h_1,1)=1$} and $c_g$ satisfies \cite{Bnt}
\begin{equation}\label{twcoc} c_g(h_1,h_2)\,c_g(h_1h_2,h_3)=c_g(h_1,h_2h_3)\,c_{g^{h_1}}(h_2,h_3)\,,\end{equation}
and therefore is a 2-cocycle for $C_G(g)$. Likewise, given a 2-cocycle $\psi$ on  $G$, define
 \begin{equation}\label{betag}\beta^\psi_g(h)=\beta_g(h):=\psi(g,h^g)\,{{\psi(h,g)^*}}\,.\end{equation}
  Then the 2-cocycle condition for $\psi$ directly yields
 \begin{align}\beta_{gg'}(h)=&\,\beta_g(h)\,\beta_{g'}(h^g)\,,\label{bichar}\\
 \beta_g(hk)=&\,\beta_g(h)\,\beta_g(k)\,\psi(h,k)\,{{\psi(h^g,k^g)^*}}\,,\label{bichar2} \end{align}
for all $g,g',h,k\in G$. In particular, for any $g\in G$, $\beta_g$ is a 1-dimensional representation of $C_G(g)$. Moreover,
when $g$ and $h$ commute, 
\begin{equation}
\beta_{g^k}(h^k)=\beta_g(h)\,.\label{conjprop}\end{equation}

As mentioned earlier,
an important class  of MTC, called the \textit{twisted group double} $\cD_\omega(G)$, is associated to $G$ and a choice of 
3-cocycle $\omega\in Z^3(G;\bbT)$. It was introduced by Dijkgraaf--Witten \cite{DW}, 
and developed in \cite{DPR}. It can be defined as the double (or centre)  of the fusion category Vec$_\omega(G)$
of $G$-graded vector spaces where the associativity constraint is defined by $\omega$, and is the category of representations of a quasi-triangular quasi-Hopf algebra \cite{Maj}. The modular data
and some initial observations about modular invariants were made in \cite{CGR}. The relation to subfactors was worked out explicitly in \cite{EP1}, and this led to the   classification  \cite{Ost} of module categories for
$\cD_\omega(G)$. This classification is rather abstract, e.g. it is very unclear which of its module categories are type 1, explicitly what are induction and restriction, what is the extended MTC, etc. As we explain in section 3.1, some clarifaction of this is made in \cite{Dav,Dav3}. Inspired by \cite{FHT}, one of the authors of this paper (DEE) established the connection of $\cD_\omega(G)$ and related structures to $K$-theory \cite{Ev}. This connection underlies this paper, and is reviewed in the next subsection.

Up to equivalence, this category $\cD_\omega(G)$ depends  only on the class $[\omega]\in H^3(G;\bbT)$ {(and of course the isomorphism class of $G$)}. Moreover, if $\alpha$ is some automorphism of $G$, then $\cD_\omega(G)$ and $\cD_{\alpha(\omega)}(G)$ are again isomorphic. The characterisation of when $\cD_\omega(G)\cong\cD_{\omega'}(G')$ is given in Corollary 1.5 of \cite{NN}.

Any $\cD_\omega(G)$ can be realised using systems of endomorphisms as follows. 
Consider any type III$_1$ factor $N$, and any subgroup of Out$(N)$ isomorphic
to $G$. Lifting that subgroup to Aut$(N)$ defines a 3-cocycle $\omega$. Then the Drinfeld or quantum double of this system  yields {the category} $\cD_\omega(G)$.
All groups $G$ and classes $[\omega]$ are realised by subgroups of Out($N$) in this way, for some $N$ --- e.g.\ $N$ can be chosen to be hyperfinite \cite{Jon1}.

The \textit{sectors} (equivalence classes of simple objects) in $\cD_\omega(G)$   are parametrised by pairs $[g,\chi]$ where $g$ is a  conjugacy class representative in $G$ and 
$\chi\in \mathrm{Irr}_{c_g}(C_G(g))$ is a projective character.
The unit  is $[1,\mathbf{1}]$. We write Ver$_\omega(G)$ for the corresponding \textit{Verlinde ring}  (Grothendieck ring) of $\cD_\omega(G)$. Ver$_\omega(G)$ is isomorphic as a ring to  the $\omega$-twisted
$G$-equivariant $K$-group ${}^\omega K^0_G(G)$, and though this is crucial to the sequel \cite{EG8} it plays no role here.

We often use in  section 4  that a group homomorphism
$\phi:G\rightarrow H$ gives rise to $K$-theoretic maps $\phi^*:\mathrm{Ver}_{\omega^\phi}(H)\rightarrow
\mathrm{Ver}_\omega(G)$ and $\phi_!:\mathrm{Ver}_\omega(G)\rightarrow \mathrm{Ver}_{\omega^\phi}(H)$, namely $[h,\tilde{\chi}]\mapsto
\sum_g[g,\tilde{\chi}\circ\phi]$ and $[g,\chi]\mapsto [\phi(g),\phi_!(\chi)]$ respectively, where the sum is over all $G$-orbit representatives $g$ with $h\in\phi^{-1}(G.g)$, and $\phi_!(\chi)$ is the wrong-way map (adjoint) of $\tilde{\chi}\circ\phi$. This is easiest to see in the bundle picture described next subsection. For example, when $\phi$ is an embedding, $\phi_!(\chi)$ is induction.

An important, though still somewhat mysterious, part of the story concerns the
modular data. The generators $S,T$ for this  SL$_2(\bbZ)$ action is, in the most general case \cite{CGR}, 
\begin{align}{S^\omega}_{[a_1,\chi_1],[a_2,\chi_2]}=&\,\frac{1}{|G|}\sum_{g_i\in \mathrm{cl}(a_i),g_1g_2=g_2g_1}{\chi_1(h_1)^*}\,{\chi_2(h_2)^*}
\,\frac{c_{g_1}(k_1^{-1},h_1)\,c_{g_2}(k_2^{-1},h_2)}{c_{g_1}(g_2,k_1^{-1})\,c_{g_2}(g_1,k_2^{-1})}\,,\label{Sfingr}\\
{T^\omega}_{[a_1,\chi_1],[a_2,\chi_2]}=&\,e^{-2\pi i c/24}\,\delta_{[a_1,\chi_1],[a_2,\chi_2]}{\chi_1(a_1)}/{\chi_1(1)}\,,\label{Tfingr}\end{align}
where $g_i=a_i^{k_i}$, and $h_1:={}^{k_1}g_2\in C_G(a_1)$, $h_2:={}^{k_2}g_1\in
C_G(a_2)$, and where $c\in 8\bbZ$ is the
central charge of the corresponding conformal net or VOA.
This SL$_2(\bbZ)$-representation  is interpreted in \cite{EG7}
using Chern characters within the bundle picture of section 2.3.

In section 4.8 we need the modular data for the cyclic group $G = \bbZ_n$. Then $H^3(\bbZ_n;\bbT) \cong \bbZ_n$, with the following explicit cocycle representatives (see e.g. \cite{HWY})
\begin{equation}\label{cyclic3cocycle}\omega_q(g_1, g_2, g_3) = \exp (2\pi i\,qg_1[(g_2 + g_3)/n]/n) \,, \end{equation}
where $q \in \bbZ_n$ parametrizes the different cohomology classes, and $[x]$ is truncation (the largest integer not greater than $x$).
The 2-cocycle $c_a^{\omega_q}(h,g)$ is  coboundary for every $q$ and $a$ (since $H^2(\bbZ_n;\bbT) = 0$), so for each $q$ there are $n^2$ sectors, which we can parametrise as
$[a,\chi_l]$ where $a\in\bbZ_n$ and $\chi_l(b) = e^{2\pi ibl/n}$ for $l \in \bbZ_n$ are the linear characters of $\bbZ_n$. The modular $T$ matrix \eqref{Tfingr} becomes
\begin{equation}
T^{\omega_q}_{[a,\chi_l],[a,\chi_l]} = e^{-2\pi i c/24}\,\exp(2\pi i\,(qa^2 + nal)/n^2) \,.\label{Tcyclic}\end{equation}

In our proof of Theorem 2 below, we use the fact that the $T$ matrix alone uniquely determines
the twist $[\omega_q]$, for $G=\bbZ_2$, $\bbZ_3$ and $\bbZ_4$. In particular, write $n_k$ ($k=0,1,\ldots,n-1$) for
the number of eigenvalues in $T^{\omega_q}$ equal to $e^{2\pi i k/n^2}$, and consider the generating polynomial
$P_{q}(x)=\sum_kn_kx^k$: then \eqref{Tcyclic} tells us

\smallskip\noindent{$\bullet$ $G=\bbZ_2$ has generating polynomials $P_0(x)=x^2+3$ and $P_1(x)=x^3+x+2$;}

\smallskip\noindent{$\bullet$ $G=\bbZ_3$ has generating polynomials $P_0(x)=    2 x^6  + 2 x^3  + 5$, $P_1(x)= 2 x^7  + 2 x^4  + 2 x + 3$, and $P_2(x)=2 x^8  + 2 x^5  + 2 x^2  + 3$;}

\smallskip\noindent$\bullet$  $G=\bbZ_4$ has generating polynomials  $P_0(x)=  2 x^{12}   + 4 x^8  + 2 x^4  + 8$, $P_1(x)= 2 x^{13}   + 2 x^{12}   + 2 x^9  + 2 x^5  + 2 x^4  + 2 x + 4$, $P_2(x)= 2 x^{14}   + 2 x^{10}   + 2 x^8  + 2 x^6  + 2 x^2  + 6$, and $P_3(x)=  2 x^{15}   + 2 x^{12}   + 2 x^{11}   + 2 x^7  + 2 x^4  + 2 x^3  + 4$.

\smallskip\noindent { For $G=\bbZ_5$ and beyond, $T^{\omega_q}$ no longer determines $\omega_q$. However, for any $G=\bbZ_n$, $T^{\omega_0}$ is the only $T^{\omega_q}$ of order exactly $n$.}

The proof of Theorem 2 also involves 3-cocycles on $G=\bbZ_2\times\bbZ_2\times\bbZ_2=:\bbZ_2^3$. The group $H^3(\bbZ_2^3;\bbT)\cong\bbZ_2^7$
is generated by 6 cocycles inflated from quotients isomorphic to either $\bbZ_2$ or $\bbZ_2^2$, together with the cocycle
\begin{equation}\omega_{iii}(a,b,c)=(-1)^{ a_1b_2c_3}\,.\label{omegaiii}\end{equation}
For a given multiplier $\beta\in Z^2(G;\bbT)$ of an abelian group $G$, the projective irreps $\chi\in
\mathrm{Irr}_\beta(G)$ all have the same dimension $d$, and there are precisely $|G|/d^2$ of them.
For $\omega_{iii}$, all $1\ne g\in\bbZ^3_2$ have nontrivial $\beta_g=\beta_g^{\omega_{iii}}\in Z^2(\bbZ^3_2;\bbT)$, each having exactly two projective irreps (both with dimension 2). For example,
the irreps $\rho_\epsilon$ with multiplier $\beta_{(1,1,1)}$ send $(1,0,0)\mapsto \left({\epsilon\atop 0}{0\atop -\epsilon}\right)$, $(0,1,0)\mapsto \left({0\atop \epsilon}{\epsilon\atop 0}\right)$, and  $(0,0,1)\mapsto \left({0\atop \epsilon i}{ -\epsilon i\atop0}\right)$, for $\epsilon=\pm 1$, and we see
that $T_{[(1,1,1),\rho_\epsilon],[(1,1,1),\rho_\epsilon]}=\epsilon i$. In fact, the MTC $\cD_{\omega_{iii}}(\bbZ_2^3)$ is braided tensor equivalent to $\cD_1(D_4)$ (see e.g. \cite{GMN}) so their $S$ and $T$ matrices coincide.

\subsection{Bundles over groupoids and twisted group doubles}

In this subsection we interpret the MTC $\cD_\omega(G)$
and its module categories using bundles over groupoids. The  classification  of module categories for $\cD_\omega(G)$ is given in Proposition 1.

A groupoid is a category whose morphisms  have both left and right inverses --- see for instance the summary in Appendix A of \cite{FHTi}. When a finite group $G$ acts on a set X, we write $X/\!/G$ for the corresponding (action) groupoid, with objects $x\in X$ and morphisms $g\in \mathrm{Hom}(x,g.x)$. A map (or homomorphism) between groupoids is a functor between the corresponding categories. An (untwisted) bundle over a groupoid is  a functor from the groupoid to the category of finite-dimensional vector spaces; we include twists, which control the projectivity of the groupoid action on the fibres. For example,  for $\psi\in Z^2(G;\bbT)$, $\psi$-twisted bundles on pt$/\!/G$ 
are projective representations of $G$ with multiplier $\psi$. 

The groupoids in this {paper} can be put in the form  $\Gamma/\!/H^L\times K^R$ where $H$ resp.{} $K$ are subgroups of a finite group $\Gamma$ which act on $\Gamma$ by left resp.{} right multiplication.
We twist by  cocycles $\widetilde{\omega}\in Z^3(\Gamma;\bbT)$, $\psi_L\in Z^2
(H;\bbT)$, $\psi_R\in Z^2(K;\bbT)$, and require the restrictions $[\widetilde{\omega}]|_H=[\widetilde{\omega}]|_K=[1]$. In this case, the twisted bundles consist of a vector space $V=\oplus_{\gamma\in\Gamma}V_\gamma$ (the total space, an $H^L\times K^R$-bimodule) carrying a projective
action of the twisted group rings $\bbC_{\psi_L}H$ and $\bbC_{\psi_R}K$, such that  $h(V_\gamma k)=V_{h\gamma k}$ as spaces.  These actions satisfy \cite{Ev} 
\begin{eqnarray}(h_1h_2)v=&\psi_L(h_1,h_2)\,{\widetilde{\omega}(h_1,h_2,{\gamma})^*}\,h_1(h_2v)\,,\label{actionleft}\\ v(k_1k_2)=&\psi_R(k_1,k_2)\,{\widetilde{\omega}({\gamma},k_1,k_2)}\,(vk_1)k_2\,,\\
h(vk)=&\widetilde{\omega}(h,{\gamma},k)\,(hv)k\,,\label{actionmixed}\end{eqnarray}
for all $h,h_i\in H$, $k,k_i\in K$, $\gamma\in\Gamma$ and $v\in V_{\gamma}$. 
A morphism  $V\rightarrow W$ between bundles is  a set of linear maps $V_\gamma\rightarrow
W_\gamma$ between the fibres which commute with the $H\times K$ action. We explain how to multiply bundles shortly.

By a slight abuse of notation, we speak of bundles over   $\Gamma/\!/\!{}_{(\tilde{\omega},\psi_L,\psi_R)}\,H^L\times K^R$ rather than $(\tilde{\omega},\psi_L,\psi_R)$-twisted bundles over $\Gamma/\!/H^L\times K^R$. We write ${}_{H}^{\,\psi_L}\cC(\Gamma)_{K}^{\,\tilde{\omega},\psi_R}$ for the category 
of bundles over the groupoid $\Gamma/\!/\!{}_{(\tilde{\omega},\psi_L,\psi_R)}\,H^L\times K^R$. When a 
twist $\psi_i$ or $\widetilde{\omega}$ is identically 1, we usually drop it. Any category of the form ${}_{H}^{\,\psi}\cC(\Gamma)_{H}^{\,\tilde{\omega},\psi}$ is fusion using the bundle product given below, provided $[\widetilde{\omega}]|_{H\times
H\times H}=1$ (see e.g. section 9.7 of \cite{EGNO}). Its double is the MTC $\cD_{\tilde{\omega}}(\Gamma)$. For $\psi=1$, and when $H$ and $\Gamma$ share no nontrivial normal subgroup, the category ${}_{H}^{\,\psi}\cC(\Gamma)_{H}^{\,\tilde{\omega},\psi}$ recovers the $A$-$A$ system of the group-subgroup subfactor $A=M\sdprod H\subset M\sdprod \Gamma=B$.

The  indecomposable bundles of $\Gamma/\!/\!{}_{(\tilde{\omega},\psi_L,\psi_R)}\,H^L\times K^R$ can be constructed  as follows.
Fix any  $H^L\times K^R$ orbit representative $\gamma\in\Gamma$ and a projective irrep
$V_\gamma$ of the stabiliser  $S=\mathrm{Stab}_{H^L\times K^R}(\gamma)=\{(h,k)\in H\times K:h\gamma =\gamma k\}\cong H\cap {}^\gamma K$ with 2-cocycle we read off from \eqref{actionleft}-\eqref{actionmixed}:
\begin{equation}c((h,k),(h',k'))=\psi_L(h,h')\,\psi_R(k'{}^{-1},k^{-1})\,{\widetilde{\omega}(hh'\gamma,k'{}^{-1},k^{-1})}\,
\widetilde{\omega}(h,h',\gamma)^*\,\widetilde{\omega}(h,h'\gamma,k'{}^{-1})^*\label{yetanothercocycle}\end{equation}
 The total space (bimodule) of the bundle associated to this pair $[\gamma,V_\gamma]$
 is the induced module Ind$_S^{H^L\times K^R}(V_\gamma)$, and $hV_\gamma k$ forms the
fibre over $h\gamma k$ for each $h\in H,k\in K$. Direct sums of bundles are defined as usual.

We multiply bundles as follows. Let $G,H,K$ be subgroups of $\Gamma$ on which {the class $[\widetilde{\omega}]$} is trivial, and choose $\psi\in Z^2(G;\bbT),\psi'\in Z^2(H;\bbT),\psi''\in
Z^2(K;\bbT)$. Consider  bundles $[a,\chi]\in {}_G^{\,\psi}\cC(\Gamma)_H^{\,\tilde{\omega},\psi'}$ and $[b,\phi]\in {}_H^{\,\psi'}\cC(\Gamma)_K^{\,\tilde{\omega},\psi''}$. 
Then $[a,\chi]\otimes[b,\phi]$ is the bundle in $ {}_G^{\,\psi}\cC(\Gamma)_K^{\,\tilde{\omega},\psi''}$
defined by
\begin{equation} \label{tensorofbundle} 
[a,\chi]\otimes[b,\phi]=\sum_{h\in G^a\cap H\backslash H/H\cap{}^bK}\left[ahb,\mathrm{Ind}_{G^{ahb}\cap H^b\cap K}^{G^{ahb}\cap K}({\chi^{hb}}{\cdot} \phi)\right]\,,\end{equation}
where $^bK=bKb^{-1}$, $G^a=a^{-1}Ga$, {$\chi^{k}(g)=\chi({}^kg)$} etc, and Ind is the induction of projective characters. In the special case where {the cocycles $\widetilde{\omega},
\psi,\psi',\psi''$ are all identically 1}, this formula reduces to  equation (2) of  \cite{KMY}; the proof of the more general case \eqref{tensorofbundle} however follows that of \cite{KMY}. 

The bundle product can be expressed in a more coordinate-free manner in terms of bimodules, as follows, using the product on the base $\Gamma$. Let $V,W$ be the
total spaces of the two bundles, i.e. the bimodules with an $\widetilde{\omega}$-projective action of the twisted group rings
$\bbC_{\psi}G^L$ and $\bbC_{\psi'}H^R$, and $\bbC_{\psi'}H^L$ and $\bbC_{\psi''}K^R$, as in \eqref{actionleft}-\eqref{actionmixed}. The total space
of the product bundle is simply the tensor product $V\otimes_{\bbC_{\psi'}H}W$ of bimodules, and the tensor product $v\otimes w$ of $v\in V_\gamma$ and $w\in W_{\gamma'}$ lies in the fibre above $\gamma\gamma'$. This bimodule picture of the bundle product  plays a crucial role in our study of 
alpha-induction in section 4.

The point is that (see e.g. \cite{Ost,Ev}) $\cD_\omega(G)$ can be identified with the category  ${}_{\Delta_G}\cC(G^2)_{\Delta_G}^{\,\tilde{\omega}}$, where 
\begin{equation}\widetilde{\omega}((g_1,g_1'),(g_2,g_2'),(g_3,g_3'))=\omega(g_1,g_2,g_3)\,\omega(g_1',g_2',g_3')^*\,\label{tildeomega}\end{equation}
{The sector $[g,\chi]$  corresponds to the bundle $[(g,1),\chi]$,}
where we identify the centraliser $C_G(g)$ with the stabiliser Stab$_{\Delta_G^{L\times R}}(g,1)$. The cocycle \eqref{yetanothercocycle} (with $\psi_L=\psi_R=1$) collapses with effort to the cocycle \eqref{2cocycle}. The braiding in this bundle picture is discussed in \cite{Ev}.
The quantum-dimension of the sector $[g,\chi]$, namely $\|\mathrm{cl}(g)\|\chi(1)$,  equals as it should the dimension of the total space divided by $|G|$.

The category-theoretic analogue of this geometric picture is module categories (see e.g. \cite{Ost1,Ost}). Identify the (fusion) category ${}_{1}\cC(\Gamma)_{1}^{\,\tilde{\omega}}$ of bundles over  the
groupoid $\Gamma/\!/\!{}_{\tilde{\omega}}\, 1$ with the fusion category Vec$_{\tilde{\omega}}(\Gamma)$ of $\Gamma$-graded vector spaces (the role of $\widetilde{\omega}$ in both categories is to define the associativity morphism for products). Let $H\le G$ be such that $[\widetilde{\omega}]|_H=[1]$, and choose any class $[\psi]\in H^2(H;\bbT)$; we can choose the cocycle $\psi$ so that $d\psi=\widetilde{\omega}$.
Then the twisted group algebra $\bbC_\psi H$ is an indecomposable algebra over Vec$_{\tilde{\omega}}(\Gamma)$; as a bundle over $\Gamma/\!/_{\tilde{\omega}}1$ it has a copy of $\bbC$ over each $h\in H$. Right $\bbC_\psi H$-modules are identified with bundles over $\Gamma/\!/\!{}_{(\tilde{\omega},\psi)}\,H^R$. The right $\bbC_\psi H$-modules form an indecomposable module
category $_1\cC(\Gamma)_H^{\,\tilde{\omega},\psi}$ over Vec$_{\tilde{\omega}}(\Gamma)$, and all module categories for Vec$_{\tilde{\omega}}(\Gamma)$ are of that form (see e.g. Example 2.1 of
\cite{Ost}). The corresponding full system $\cC_{\cM}^*$ for the module category $\cM={}_1\cC(\Gamma)_H^{\tilde{\omega},\psi}$ is $_H^{\,\psi}\cC(\Gamma)_H^{\,\tilde{\omega},\psi}$.

Continuing this example, we identify  $\bbC_\psi H$-bimodules {in the category Vec$_{\tilde{\omega}}(\Gamma)$} with bundles over $\Gamma/\!/\!{}_{{(\tilde{\omega},\psi,\psi)}}\,H^L\times H^R$, i.e. we identify the fusion categories Bimod$_{\mathrm{Vec}_{\tilde{\omega}}(\Gamma)}(\bbC_\psi H)$ and ${}_{H}^{\,\psi}\cC(\Gamma)_{H}^{\,\tilde{\omega},\psi}$. Theorem 3.1 of \cite{Ost} says the indecomposable module categories over ${}_{H}^{\,\psi}\cC(\Gamma)_{H}^{\,\tilde{\omega},\psi}$ are ${}_{H}^{\,\psi}\cC(\Gamma)_{H_1}^{\,\tilde{\omega},\psi_1}$, 
 where $H_1\le\Gamma$, $[\widetilde{\omega}]|_{H_1}=[1]\in H^3(H_1;\bbT)$, and $\psi_1\in Z^2(H_1;\bbT)$. Their bundles are $(\bbC_\psi H,
\bbC_{\psi_1}H_1)$-bimodules. Two pairs $(H_i,\psi_i)$ yield equivalent module categories if they are conjugate, i.e. if there is
some $g\in\Gamma$ such that $H_1=H_2^g$ and $[\psi_2(h,h')]=[\psi_1(h^g,h^{\prime g})]$.

This leads to  Ostrik's characterisation of the module categories of $\cD_\omega(G)$:

\medskip\noindent{\textbf{Proposition 1.} (\cite{Ost}, Theorem 3.2) \textit{The indecomposable module categories of $\cD_\omega(G)$ are, up to equivalence,  ${}_{\Delta_G}\cC(G^2)_H^{\,\tilde{\omega},\psi}$, where $H\le G\times G$, $[\tilde{\omega}]|_H=[1]$,
and $[\psi]\in H^2(H;\bbT)$. The module categories ${}_{\Delta_G}\cC(G^2)_H^{\,\tilde{\omega},\psi}$ and ${}_{\Delta_G}\cC(G^2)_{H'}^{\,\tilde{\omega},\psi'}$ are equivalent iff $(H,[\psi])$
and $(H',[\psi'])$ are conjugate.}}\medskip

 For example, the identity modular invariant $\cZ=I$ corresponds to the choice
$H=\Delta_G$ and $\psi=1$. {More generally,  for $H=\Delta_G$,} $\psi$ plays the role of discrete torsion, twisting $I$
by the associated cocycles $\beta_g^\psi$, {but for more general $H$,  $\psi$ can act in a  much
more complicated way}, and the general relation between pairs $(H,[\psi])$ and  modular invariants $\cZ$ is subtle, as we see in Theorem 1.

We can describe the $\cD_\omega(G)$-module category 
 ${}_{\Delta_G}\cC(G^2)^{\,\tilde{\omega},\psi}_H$ corresponding to $(H,\psi)$ in more detail as follows.
Let $(h_+,h_-)\in H$ act on $g\in G$ on the right by $g\mapsto h_-^{-1}gh_+$.
The indecomposable bundles over $G^2/\!/\!_{(\omega,1,\psi)} \,\Delta_G^L\times H^R$ are parametrised by
 pairs $[g,\chi]$ for $g$  a representative of an $H$-orbit 
in $G$, and $\chi$ a projective character  of the stabiliser $\{(h_+,h_-)\in H:h_-g=gh_+\}$ in $H$ of 
$g$. The special bundle, corresponding in the subfactor
language to the inclusion $\iota:N\subset M$ of factors, is
$\iota=[1,\mathbf{1}]$ while its conjugate $ \overline{\iota}$ is the bundle $[1,\mathbf{1}']$ over $G^2/\!/\!_{(\omega,\psi,1)} \,H^L\times \Delta_G^R$, where $\mathbf{1},\mathbf{1'}$ are certain 
projective {characters} on $H\cap \Delta_G$ ($c^\omega_1$ equals 1, but $\psi$ in
general contributes a nontrivial 2-cocycle). Then $\overline{\iota}\iota$ is the algebra $A$ for the module category and the dual canonical endomorphism $\theta$ for the subfactor, and
can be computed using \eqref{tensorofbundle} to be
\begin{equation}\label{theta}\theta=A=\sum_{h\in\Delta\cap H\backslash H/(\Delta\cap H)}
\left[\Delta h\Delta,\mathrm{Ind}_{\Delta^h\cap H\cap \Delta}^{\Delta^h\cap\Delta}\mathbf{1}\otimes\mathbf{1}'\right]\in{}_{\Delta}\cC(G^2)^{\,\omega}_\Delta\,,\end{equation}
where we write $\Delta$ for $\Delta_G$, as in \cite{Ev} in the untwisted case. 
The  corresponding full system  $\MXM$ is  ${}_H^{\,\psi}\cC(G^2)_H^{\,\tilde{\omega},\psi}$ \cite{Ev}.  Alpha-induction, sigma-restriction, the product in $\MXM$, the Ver$(G)$-module structure
of $\MXN$, etc are nontrivially modified {by $\psi$}.

\section{Results and examples}

\subsection{{The type 1 module categories of twisted group doubles}}

Let $G$ be finite {and $\omega\in Z^3(G;\bbT)$. Recall the classification (Proposition 1)
of module categories for the twisted group double $\cD_\omega(G)$ in terms of pairs $(H,\psi)$.} In this subsection we give {our first result,} Theorem 1, which is our answer to  the following question of  Ostrik \cite{Ost}:

\medskip\noindent\textbf{Problem.} \textit{Which pairs {$(H,\psi)$ in Proposition 1 correspond to type 1 (i.e. extension-type)
module categories}?} \medskip

Ostrik also asked to determine explicitly the modular invariant $\cZ$ associated to
each pair $(H,\psi)$ (not just the type 1 ones). We address that in \cite{EG8}, which gives a remarkably simple and conceptual $KK$-theory description of the module categories for finite group doubles. As we show in \cite{EG7}, this $KK$-theory description extends to the MTC associated to lattices and loop groups.

Alpha-induction and the modular invariant in the case where $\omega=1$ and $\psi = 1$, with $\Delta \le H$, was handled in
\cite{Ev}, and the modular invariant, full system, nimrep etc for any $H$ where $\omega =1$ and $\psi = 1$ was handled in \cite{EG3}. For arbitrary $\omega$, Davydov--Simmons \cite{Dav3} in Theorem 2.19 identify which pairs $(H,\psi)$ are type 1 (for $\omega=1$, this was done in Theorem 3.5.1 of   \cite{Dav}), and the MTC of local modules is given in Theorem 2.20; Theorem 4.8 gives the modular
invariant matrix, though in the basis of commuting pairs $(g,h)\in G^2$ rather than the natural basis of sectors $[g,\chi]$ (for $\omega=1$, this was done in Theorem 3.5.3 of \cite{Dav}). (This commuting pairs basis simplifies many calculations because the modular data \eqref{Sfingr},\eqref{Tfingr} become generalised permutation matrices \cite{Bnt}, but e.g. integrality is lost.)
 To our knowledge, none of alpha-induction, the full system, nor sigma-restriction is given in \cite{Dav3,Dav}.

For type 1 module categories,  sigma-restriction of local modules determines the modular invariant (but not conversely). The expression in \cite{Dav3} for the modular
 invariant associated to  $(H,\psi)$ takes the complicated  form
\begin{equation}\cZ_{[f,f'],[g,g']}=\sum_{(y,y')\in Y}\frac{c^f(g^y,y^{-1})\,c^{f'}(g^{\prime y'},y^{\prime-1})^*}{c^f(y^{-1},g)\,c^{f'}(y^{\prime-1},g')^*}\frac{\psi((f^y,f^{\prime y'}),(g^y,g^{\prime y'}))}{\psi((g^y,g^{\prime y'}),(f^y,f^{\prime y'}))}\,,\label{DSmodinv}\end{equation}
 for $c^f(a,b)=\omega(a,b,f)^*\omega(a,{}^bf,b)\omega(f^{ab},a,b)$, where $Y$ is a set of coset representatives for $\{(y,y')\in G^2\,|\,(f^y,f^{\prime\,y'})\in H\}/H$. Again, this is expressed in the wrong basis, which is why the matrix entries aren't in $\bbZ_{\ge 0}$. Compare this with the expression $\cZ=bb^t$, where $b$ is the matrix for the linear map Res$:Ver_{\bar{\omega}\bar{\omega}^\psi}(\overline{K})\rightarrow \mathrm{Ver}_\omega(G)$ given by \eqref{sigresgen} and expressed in the natural bases $[\overline{k},\overline{\chi}]$ and $[g,\chi]$. 
Again, the more fundamental quantity is sigma-restriction $b$, but that is obscure in \eqref{DSmodinv}.The equation $\cZ=bb^t$ applies to type 1 $(H,\psi)$; we give an
analogous expression in \cite{EG8}  for the modular invariant for general $(H,\psi)$, again using only inductions of characters and decomposition of conjugacy classes of a group into conjugacy classes of a subgroup.

\medskip\noindent\textbf{Theorem 1.} \textit{Let $G$ be any finite group and $\omega\in Z^3(G;\bbT)$ be any 3-cocycle, and define $\widetilde{\omega}$ by \eqref{tildeomega}.}

\medskip\noindent\textbf{(a)} \textit{Let $H\subseteq G^2$ satisfy $[\widetilde{\omega}]|_H=[1]$, and choose $\psi\in Z^2(H;\bbT)$. Then $(H,\psi)$  is a type 1 pair for $\cD_\omega(G)$ iff (up to conjugation in $G^2$) $H=\Delta_K(1\times N)$ for some subgroup $K$ of $G$ and normal subgroup  $N$ of $K$, with $[\psi]|_{\Delta_K}=[1]$ and $\beta^\psi_{N\times 1}(1\times N)=1$. The condition on $\psi$ is equivalent (up to adjustment by coboundary) to requiring that $\psi$  satisfy both $\psi(\Delta_K,H)=1$
and $\psi((n,1),(1,n'))=\psi((1,n'),(n,1))$ for all $n,n'\in N$.}

 \medskip\noindent\textbf{(b)}  \textit{Choose any type 1 pair $(H,\psi)$, with $\psi$ as in the last sentence of (a).}
 
\begin{itemize}\item[{(i)}] \textit{The associated algebra $A$ is described below. The category of right $A$-modules (the module category) is ${}_{\Delta_G}\cC(G^2)^{\tilde{\omega},\psi}_{H}$, and the category of $A$-$A$-bimodules (the full system) is ${}_{H}^{\,\psi}\cC(G^2)^{\tilde{\omega},\psi}_{H}$. Alpha-induction $\alpha:\cD_\omega (G)\rightarrow{}_{\Delta_G}\cC(G^2)^{\tilde{\omega},\psi}_{H}$ is
 \begin{equation}\alpha^H([g,\chi])=\sum_{k\in C_G(g)\backslash G/K}\bigl[g^k, \mathrm{Ind}_{ C_K(g^k)}^{C_H(g^k)}\,\chi^{k}\bigr]
 \,,\label{alphagen}\end{equation}
where $C_H(g)=\{k\in K\,|\,g^k\in N\}$. Write \eqref{alphagen} symbolically as $\sum [g',\widetilde{\chi}]$.
The alpha-inductions $\alpha_\pm:\cD_\omega (G)\rightarrow{}_{H}^{\,\psi}\cC(G^2)^{\,\tilde{\omega},\psi}_{H}$ are
 $\alpha_+[g,{\chi}]= \sum [(g',1),\widetilde{\chi}_+]\,$ and
 $\alpha_-[g,{\chi}]=\sum [(1,g^{\prime-1}),\widetilde{\chi}_-]\,,$
 where we 
 write $\widetilde{\chi}_\pm(h_+,h_-)=\widetilde{\chi}(h_\pm,h_\pm)$.}
 
\item[{(ii)}] \textit{Sigma-restriction ${}_{\Delta_G}\cC(G^2)_H^{\tilde{\omega},\psi}\rightarrow {}_{\Delta_G}\cC(G^2)^{\,\tilde{\omega}}_{\Delta_G}$ sends the bundle $[(g,1),\chi]$ to
$\sum_n[(ng,1),\chi']$, where the sum is as in the decomposition $\mathrm{cl}_K(g)N=\cup_n\mathrm{cl}_K({ng})$ into  orbits $\mathrm{cl}_K(g')$  of $K^{adj}$ in $G$, and $\chi'$ is defined in \eqref{rho'}. The  indecomposable local bundles in ${}_{\Delta_G}\cC(G^2)_H^{\tilde{\omega},\psi}$ are those bundles $[(k,1),{\chi}]$ for which $k\in K$ and $\chi$ satisfies both $|\chi(k,k)|=|\chi(1,1)|$ and
\eqref{local} for all $n\in N$.
}

\item[{(iii)}] \textit{
The modular tensor category of local modules for $(H,\psi)$  is braid-equivalent to the twisted group double $\cD_{\bar{\omega}\bar{\omega}^\psi}(\overline{K})$, where $\overline{\omega}$ and $\overline{\omega}^\psi$ are 3-cocycles 
on $\overline{K}:=K/N$ defined in Proposition 2 (with $G$ there replaced with $K$). Sigma-restriction
Ver$_{\bar{\omega}\bar{\omega}^\psi}(\overline{K})\rightarrow\mathrm{Ver}_\omega(G)$
is $\iota_!\circ\pi^*$ where $\pi^*$ resp.\ $\iota_!$ are the natural $K$-theoretic maps coming from the quotient $\pi$ and inclusion $\iota$ respectively.}
\end{itemize}

$K$-theoretic maps like $\pi^*:\mathrm{Ver}_{\bar{\omega}\bar{\omega}^\psi}(\overline{K})\rightarrow\mathrm{Ver}_\omega (K)$ and $\iota_!:\mathrm{Ver}_\omega (K)\rightarrow \mathrm{Ver}_\omega (G)$ in Theorem 1(b)(iii) are discussed in section 2.2. Explicitly, choosing coset representatives $t_{\bar{k}}\in\overline{k}\in K/N$, the sigma-restriction $\iota_!\circ\pi^*$ in Theorem 1(b)(iii)  becomes
\begin{equation}\label{sigresgen}
[\overline{k},\overline{\chi}]\mapsto\sum_k[k,\mathrm{Ind}_{C_K(k)}^{C_G(k)}\tilde{\phi}\,\overline{\chi}\circ\pi]\end{equation}
where  the sum over $k$ is as in $\cup_k\mathrm{cl}(k)=\pi^{-1}(\mathrm{cl}(\overline{k}))$, and $\widetilde{\phi}(h)=\beta_{(k,k)}(ht_{hN}{}^{-1},1)^*$. 

 Ostrik's guess \cite{Ost} for the type 1 module categories missed the $\psi\ne 1$ possibilities. The smallest $G$ for which these arise are the dihedral group $D_4$ and the alternating group $A_4$ --- see section 3.3 for an example.
We do not know whether inequivalent module categories $(H,\psi)$ necessarily have
distinct modular invariants $\cZ$ --- this would have seemed very unlikely except that it is true when $\cZ$
is a permutation matrix \cite{EG8}.

The MTC of local modules goes by other names in the literature: the ambichiral or neutral system, and  more generally the type 1 parents, in subfactor theory; the left and right centres in categorical treatments; the extended system in CFT or VOA. The phrase we use should make sense in all frameworks. 

Of the two contributions $\overline{\omega}$ (which depends on $\omega$) and $\overline{\omega}^\psi$ (the part depending on $\psi$) to the 3-cocycle twist in
Theorem 1(b)(iii), the more important for us is $\overline{\omega}^\psi$. Our proof of reconstruction (Theorem 3)
requires delicate control over $\overline{\omega}^\psi$ --- in fact we get reconstruction (Theorem 3) for arbitrary $\omega$,
even if Theorem 1 were to hold  only for trivial twist $\omega=1$. A crucial part of our proof connects $\overline{\omega}^\psi$ to the thesis of Vaughan Jones \cite{Jon2}.
A $Q$-kernel is a homomorphism  from a group $Q$ to the outer automorphism group Out$(M) = \mathrm{Aut} (M)/\mathrm{Int} (M)$  of a factor $M$. Any $Q$-kernel  gives rise to a 3-cocycle on  $Q$. If one has an action of $K$ by automorphisms on  a factor, which is inner when restricted to subgroup $N$,
then one has a $K/N$-kernel and hence  a 3-cocycle on $K/N$.
This is what Jones exploits and computes explicitly. The structure coming along with this framework is necessary for our proof of Theorem 3.

We find in section 4.3 that the algebra $A$ associated to the type 1 module category
$(H,\psi)$ of Theorem 1 consists of a copy $A_{(gn,g)}$ of $\beta^\psi_{(n,1)}\otimes \cF(C_G(n)/C_K(n))$ at each point $(gn,g)\in \Delta_G(1\times N^G)$,  
where $N^G$ is the set $\cup_{g\in G}N^g$ and $\cF(G'/K')$ denotes the algebra of functions on the set $G'/K'$ of cosets. Multiplication in $\cF(G'/K')$ is
pointwise, and as a $G'$-representation is equivalent to
Ind$_{K'}^{G'}\mathbf{1}$. The group $\Delta_G$ acts on the left and right of $A$ as in \eqref{actionleft}-\eqref{actionmixed}, making $A$ into a bundle in ${}_{\Delta_G}\cC(G^2)^{\,\tilde{\omega}}_{\Delta_G}$. The multiplication $\mu:A_{(1,1)}\otimes A_{(n,1)}\rightarrow A_{(n,1)}$ is
that of $\cF(G/K)$ restricted to $C_G(n)$, and is extended to $\mu:\bbC_{(gn,g)}\otimes\bbC_{(g'n',g')}\rightarrow\bbC_{(gg'n^{g'}n',gg')}$ by equivariance 
\begin{equation}\label{multequiv}g\,\mu(v\otimes w)\,h=\mu(gv\otimes wh)\,.\end{equation}

It is elementary to specialise Theorem 1 to  cyclic groups.

\medskip\noindent\textbf{Corollary 1.} \textit{Let $G\cong \bbZ_n$ and write $\omega=\omega_q$ as in \eqref{cyclic3cocycle}. 
Then the  type 1 module categories correspond to any choice of integers $m,m',a\in\bbZ_{>0}$, where $m'|m|n$, $a\le m'$,
$m'|q$ and $(m')^2|am$. The MTC of local bundles is a twisted double of $\bbZ_{m/m'}$.} \medskip

Here, if we write $G=\langle g\rangle$, then $K=\langle g^{n/m}\rangle\cong \bbZ_{m}$ and $N=\langle
g^{n/m'}\rangle\cong \bbZ_{m'}$. The 2-cocycle $\psi$ comes from $a$ and is given explicitly in
section 4.6, along with the proof. See also the discussion of \cite{EMS} in section 3.2. Cyclic $G$ is
so easy because all sectors in $\cD_\omega(G)$ are invertible.

The maximal extensions of $\cD_\omega(G)$ occur when the MTC of local modules is trivial, i.e. is Vec$_\bbC$. Such type 1 module categories are said to \textit{trivialise} $\cD_\omega(G)$, for this reason. They are classified by an analogue of Proposition 1 \cite{Ost}, but again the branching rules etc  are obscure. In section 4.7 we recover this classification as a corollary of Theorem 1, and with it the branching rules etc.

\medskip\noindent\textbf{Corollary 2.} \textit{Choose any finite group $G$ and any 3-cocycle $\omega\in Z^3(G;\bbT)$. The type 1 module categories  of $\cD_\omega(G)$ 
whose MTC of local modules is $\mathrm{Vec}_\bbC$, correspond to any subgroup  $K\le G$ and  $\psi'\in Z^2(K;\bbT)$ such that ${\omega}|_{K}$  is coboundary. Then
\begin{equation}1\mapsto\sum_k\,\bigl[k,\beta^{\psi'}_k\cdot\mathrm{Res}_{C_G(k)}\mathrm{Ind}_K^G\mathbf{1}\bigr]\,,\end{equation}
is  the associated branching rule, where the sum is over conjugacy class representatives of $K$ by $G$, and $\beta^{\psi'}$ is defined in \eqref{betag}. In the notation of Theorem 1(a), this module category
has $H=K\times K$ and cocycle $\psi\in Z^2(K\times K;\bbT)$  given by} 
\begin{equation}\psi((k_1,k_2),(k_1',k_2'))=\psi'(k_2'{}^{-1},k_2^{-1}k_1)\,\psi'((k_2^{-1}k_1)^{k_2'},k_2'){}^*\,\psi'((k_2^{-1}k_1)^{k_2'},k_2'{}^{-1}k_1')\,.\end{equation}

In the physics literature, two of these maximal extensions are
singled out: the pair $(K,\psi')=(1,1)$ (the original holomorphic theory) and the pair $(K,\psi')=(G,1)$ (what is unfortunately called there the orbifold of the original theory by $G$). The physics literature realises that
the latter doesn't always exist --- in particular so-called \textit{level-matching} must be satisfied. This (necessary but not sufficient) condition requires that the resulting modular invariant be invariant under $T$. We see no reason to restrict to those two maximal extensions. Moreover level-matching must be replaced with the $[{\omega}]|_K=[1]$ condition, which implies it; in general it is only  for cyclic $K$
that level-matching is equivalent to that  condition  on $\widetilde{\omega}$.

For example, Corollary 2 for $G=\bbZ_n$ reduces to $m'=m$ dividing $q$, and $a=0$.

{Although Proposition 1 is abstract nonsense (for this reason Ostrik asked his questions), it is significant that the formulae of Theorem 1} 
yield explicit expressions. Indeed, all alpha-inductions, sigma-restrictions and modular invariant entries are computed by nothing more complicated than induction and restriction of (possibly projective) characters, and decomposition of conjugacy classes of a group into conjugacy classes of a subgroup. All coefficients are manifestly in $\bbZ_{\ge 0}$. Before turning to the proofs in section 4, we discuss  some applications of Theorem 1 in section 3.2 and collect several concrete examples  in section 3.3.

\subsection{{VOAs, conformal nets and finite group doubles}} 

This subsection describes some deep consequences of Theorem 1  for conformal nets of factors. Similar comments apply to  vertex operator algebras (VOAs), provided the following conjecture holds. Recall that a holomorphic conformal net $\cA$ (resp. VOA $\cV$) is a completely rational net resp. VOA with trivial representation theory: Rep$(\cA)\cong \mathrm{Vec}_{\bbC}\cong\mathrm{Mod}(\cV)$. When $G$ is a finite group of automorphisms, the space $\cV^G$ of fixed points is automatically a VOA, called the \textit{orbifold} of $\cV$ by $G$.

\medskip\noindent\textbf{Conjecture 1.} \textit{Let $\cV$ be a holomorphic VOA and $G$ be a finite group of
automorphisms of $G$. Then the orbifold $\cV^G$ is also completely rational.}\medskip

In fact it is generally believed that $\cV^G$ is completely rational whenever $\cV$ is completely rational (and not necessarily holomorphic) and $G$ is finite. This has been proved \cite{Miy,CM} when $G$ is solvable, and the analogue is known to be true for completely rational conformal nets \cite{Xu}. 

\medskip\noindent\textbf{Conjecture 2.} \textit{Let $\cV$ be a holomorphic VOA and $G$ a finite group of automorphisms. Assume $\cV^G$ is completely rational (Conjecture 1). Then
the MTC Rep($\cV^G)$ is $\cD_\omega(G)$ for some
3-cocycle $\omega\in Z^3(G;\bbT)$.}\medskip

This also is generally believed. Proposition 5.6 of \cite{DRX} establishes this for all cyclic $G$, using
intrinsically VOA methods. The main theorem of \cite{Kir} shows this for any $G$ but trivial
twist $\omega$  (i.e. all 2-cocycles arising are 1) --- his proof is categorical, so works equally
for conformal nets and (assuming Conjecture 1) VOAs. For conformal nets, Conjecture
2 is proven for all $G,\omega$ (Theorem 5.14 in \cite{Mug}). We will see shortly that the converse also holds: any completely rational conformal net with Rep$(\cA)\cong\cD_\omega(G)$ can be obtained as a $G$-orbifold of a holomorphic conformal net.

For the remainder of this subsection unless otherwise stated, whenever we discuss VOAs, we will assume both Conjecture 1 and 2 hold; assuming these (and often just Conjecture 1), all statements made for conformal nets also hold for
VOAs.

As explained in the Introduction, a natural question is to describe all conformal (a.k.a.{} local) extensions of the orbifold
$\cA^G$. Such questions  seem very difficult to answer at present, unless one uses categorical methods:
 the conformal extensions of $\cA^G$ (or $\cV^G$ provided it is completely rational) correspond to the type 1 module categories \cite{Kir,HKL,CKM} in Rep$(\cA^G)\cong\cD_\omega(G)$. Thus Theorem 1 gives all such extensions $(\cA^G)^e$ of $\cA^G$, identifies each MTC Rep$((\cA^G)^e)$ as a twisted group double, 
and expresses the decomposition (branching rules) of $(\cA^G)^e$ and its
representations $\pi\in\mathrm{Rep}((\cA^G)^e)$ into a direct sum of $\cA^G$-representations. Equivalent module categories ${}_{\Delta_G}\cC(G^2)_H^{\tilde{\omega},\psi}$ (i.e. conjugate pairs $(H,\psi)$) correspond to extensions of $\cA^G$ equivalent as $\cA^G$-modules (and also as conformal nets),
 but inequivalent pairs $(H,\psi)$ may correspond to equivalent conformal nets (or VOAs) --- see e.g. the $\cV(A_2\oplus E_6)$ lattice VOA example given next subsection. 

Some of these extensions are clear. The \textit{quantum Galois correspondence}  \cite{HMT,DoMa} says that
all conformal nets $\cB$ between $\cA^G$ and $\cA$ are in natural bijection with subgroups $K$
of $G$,  namely $K\leftrightarrow \cA^K$. These correspond to the 
type $1$ module categories with $(H,\psi)=(\Delta_K,1)$. But Theorem 1 tells us that there are several other conformal extensions $(\cA^G)^e$ of $\cA^G$ which aren't subnets of $\cA$. Much more challenging is
to obtain these from traditional conformal net or VOA constructions. The plethora of extensions contained in Theorem 1 provides a graphic
example of what seems to be a serious {general} challenge for VOAs and conformal nets: to find new construction methods.
See e.g. \cite{EG2} for further VOA candidates  (realised by braided subfactors and MTC) that have not yet been constructed.

The simplest class of these less obvious extensions is as follows. Let $\cA$ be holomorphic, $G$  a finite group of automorphisms and Rep$(\cA^G)\cong\cD_\omega(G)$,  as before.
Assume for simplicity here that $\omega=1$ --- nontrivial $\omega$ is handled by Theorem 1, and excludes some of the following extensions. There exists
a completely rational conformal net $\cA^G_N$ containing $\cA^G$, for each normal subgroup $N\triangleleft G$.
In this case the category of $\cA^G_N$-modules is the double $\cD_{1}({G}/N)$. Any simple $\cA^G_N$-module
decomposes as a $\cA^G$-module into the direct sum \eqref{sigma2i}. More generally, one can twist  these  
$\cA^G_N$ by Jones' characteristic invariant (the $\psi$ in Theorem 1(a)), and this can introduce a nontrivial 3-cocycle $\overline{\omega}^\psi$ even when one wasn't there before. The general extension of $\cA^G$ is a combination of the $\cA^H$ and these $\cA^G_{N,\psi}$.
  
  Corollary 2  gives all holomorphic (i.e. maximal) extensions of $\cA^G$ (or $\cV^G$). There may
  only be one of these, namely  the original  $\cA$ or $\cV$ itself, corresponding to $K=1$. Conversely, any conformal net $\cA$ with Rep$(\cA)\cong
\cD_\omega(G)$, is an orbifold by $G$ of a holomorphic net: choose the holomorphic extension $\cA^{hol}$
of $\cA$ associated to $K=1$ and $\psi= 1$ in the notation of Corollary 2; then as a sum of sectors $\cA^{hol}=\oplus_{\pi\in\mathrm{Irr}(G)}
\pi\otimes [1,\pi]$ using obvious notation; then each summand $\pi\otimes[1,\pi]$ manifestly carries a $G\times \cA$-action,
and it is evident that $(\cA^{hol})^G=[1,\mathbf{1}]=\cA$.

The case of cyclic $G$ was studied recently in  \cite{EMS}. This paper plays an important role both in the classification of
holomorphic VOAs with central charge 24 (they construct 5 new ones), and in the proof of Generalised Moonshine. As we know,
Conjectures 1 and 2 are known to hold for $G$ cyclic. Their Theorem 5.16 describes all type 1 module categories --- because all sectors are simple-currents here it is straightforward to get
the complete story.
Holomorphic extensions of cyclic orbifolds, which share the additional key property of the Monstrous Moonshine VOA that there are few states with small conformal weights, have recently been studied in \cite{GK}.

The best studied examples of orbifolds are \textit{permutation orbifolds}. Let $\cA$ be any conformal net, and $\cV$ be any VOA. Let $G$ be any subgroup of the symmetric group $S_k$. Then $(\cA^{\otimes k})^G$ and $(\cV^{\otimes k})^G$ are called permutation orbifolds --- here $G$ acts by $v_1\otimes\cdots \otimes v_k\mapsto v_{\pi 1}\otimes\cdots \otimes v_{\pi k}$ in the $k$-fold tensor product. 

We are interested in the special case where $\cA$ is holomorphic.  In this case the central charge $c$ necessarily lies in $8\bbZ_{\ge 0}$. 
The categories Rep$((\cA^{\otimes k})^G)$ and (assuming Conjectures 1 and 2) Mod$((\cV^{\otimes k})^G)$
is $\cD_\omega(G)$ up to some twist $\omega\in Z^3(G;\bbT)$. In Theorem 2 we identify
this twist $\omega$ for any such permutation orbifold.

\medskip\noindent\textbf{Theorem 2.} \textit{Let $\cA$ be a holomorphic conformal net of central charge $c$, and let $G$ be a subgroup of some symmetric group $S_k$. Let $(\cA^{\otimes k})^G$ be the corresponding permutation orbifold, and let $\cD_\omega(G)$ be its category of modules. Then $[\omega]$ is the restriction to $G$ of $[\omega^{(3)}_{c\,\, (\mathrm{mod}\,\, 3)}]\in H^3(S_k;\bbT)$, where the 3-torsion cocycle $\omega^{(3)}_q$ is defined in section 4.8. In particular, if  24 divides $c$ or 3 does not divide the order $|G|$ of $G$, then the twist $\omega$ is trivial.}\medskip

We thank Marcel Bischoff for correspondence on the possibility that $[\omega]\ne [1]$. The expectation had been that $[\omega]=1$ always \cite{Mug1,Dav2}. Both Bischoff and  Johnson-Freyd (the latter in \cite{JF})  also showed the 3-cocycle could be nontrivial in general. Example 2.1.1 in \cite{JF} announces the general result, corresponding to our Theorem 2. A conjecture-free generalisation to VOAs is given in Theorem 4(a) below.

The nontrivial twist can be detected by the character vectors.
For example, fix some holomorphic theory $\cA$ at central charge $c\in 8\bbZ$ and some $n\ge 1$, with graded dimension $\chi_\cA(\tau)$, and write $\pi=(12\cdots n)\in S_n$. Consider the $\bbZ_n$-permutation orbifold $\cA_\pi=(\cA^{\otimes n})^{\langle \pi\rangle}$. 
Then $\cA_\pi$ has graded dimension 
\begin{equation}\label{permorbchar}\chi_{\cA_\pi}(\tau)=\sum_{d|n}\frac{\phi(d)}{n}\,\chi_\cA(d\tau)^{n/d}\,\end{equation}
where $\phi(d)$ is the Euler totient and the sum is over all divisors of $n$ (see equation (4.15) of \cite{BHS} for an analogous expression for the character of any $\cA_\pi$-module). Now, $\chi_\cA(-1/\tau)$ must be a linear combination with nonzero coefficients of the characters of $\cA$-representations, and so from it we can read off all values appearing in the modular $T$ matrix for $\cA$. We see that 
$$\chi_{\cA_\pi}(-1/\tau)=\sum_{d|n}\frac{\phi(d)}{n}\,\chi_\cA(\tau/d)^{n/d}\in \sum_{d|n}q^{-nc/24d^2}\bbC[[q^{1/d}]]\,.$$
The central charge of $\cA_\pi$ is $nc$, and we see that $e^{2\pi i nc/24}T$ for $\cA_\pi$ has order $n$, unless 
both $3|n$ and $3$ doesn't divide $c$, in which case the order is $3n$. But $T$ should match $T^{\omega_q}$ in \eqref{Tcyclic} for some $0\le q<n$, so $q=0$ unless both $3|n$ and $3\!{\not{}|c}$,
in which case $q=cn/3$. 

Perhaps the biggest question on the interface of conformal nets/VOAs and categories is \textit{reconstruction}: is every MTC
realised by a completely rational VOA or  conformal net? The first nontrivial place to look, presumably, is to the finite group categories 
$\cD_\omega(G)$, for arbitrary $G$ and $\omega$. These would correspond to orbifolds by
$G$ of a holomorphic VOA or net. 

It had been expected for some time that permutation orbifolds would realise
any $\cD_1(G)$, and this is now proven (Theorem 2).
But until now, it has been far from clear that $\cD_\omega(G)$ for all 3-cocycles $\omega$ and all finite groups $G$,
can likewise be realised as a representation category. Indeed, the 3-cocycles appearing in Moonshine
all seem to have order dividing 24. 

\medskip\noindent\textbf{Theorem 3.} \textit{For any finite group $G$ and any 3-cocycle
$\omega\in Z^3({G};\bbT)$, there exists a completely rational conformal net
$\cA$ whose representation  category Rep$(\cA)$  is $\cD_\omega(G)$.}\medskip

For VOAs,  weaker statements hold (we assume neither Conjecture 1 nor 2 for the proof):

\medskip\noindent\textbf{Theorem 4(a)} \textit{Let $\cV$ be a holomorphic VOA of central charge $c$ divisible by 24, and let $G$ be a solvable subgroup of some symmetric group $S_k$. Let $(\cV^{\otimes k})^G$ be the corresponding permutation orbifold. Then its category of modules is $\cD_1(G)$.}

\medskip\noindent\textbf{(b)} \textit{When $G$
is solvable, then for any $\omega\in Z^3(G;\bbT)$ there exists a completely rational VOA $\cV$
with $\mathrm{Mod}(\cV)\cong \cD_\omega(G)$.}\medskip

When $G$
is not solvable, a completely rational VOA  realising $\cD_\omega(G)$ for any $\omega\in Z^3(G;\bbT)$
exists provided  we know $\cV^G$ is completely rational for some holomorphic VOA $\cV$. We prove Theorems 2, 3 and 4 in section 4.

\subsection{{Examples}}

Consider any type 1 module category, i.e. any pair $(H,\psi)$  in Theorem 1(a). The simple objects in the extended theory, thanks to Theorem 1(b)(iii), are
 $[\overline{k},\overline{\chi}]$, where $\overline{k}\in\overline{K}=K/N$, and $\overline{\chi}$ is a projective character of $C_{\overline{K}}
(\overline{k})$ with 2-cocycle  \eqref{2cocycle}. Sigma-restriction, and with it the modular invariant, are determined from \eqref{sigresgen}.

Let's turn next to some concrete examples. The trivial module category  for $\cD_\omega(G)$, corresponding 
to the unextended conformal net or VOA, and the diagonal modular invariant $\cZ=I$,  is associated to the choice $H=\Delta_G$ and $\psi=1$.

Another example which works for any $\cD_\omega(G)$ is  $H=1\times 1$ and $\psi=1$. In this case, there
is only one local module (so the extended theory is holomorphic), with sigma-restriction 
$[1,\mathbf{1}]\mapsto \sum_\pi \mathrm{dim}\,\pi\,[1,\pi]\in\mathrm{Ver}_\omega(G)$, where the sum is over all (linear)
characters $\pi\in\mathrm{Irr}(G)$. Note that the 2-cocycle multiplier for these $\pi$ is trivial
even when $\omega$ is nontrivial. This example appeared in equation (1.19) of \cite{EG2}.

Another holomorphic example (which works though only for $\omega=1$) is $H=G\times G$. Take $\psi=1$  for simplicity. Then  sigma-restriction sends the unit to $\sum_g[g,\mathbf{1}]$ where the sum is over conjugacy class
 representatives $g$. This example appeared in equation (1.20) of \cite{EG2}.
 
For comparison, a simple example which is not type 1 is $H=G\times 1$ and $\psi=1$. Then we have two type 1 parents (a.k.a. left and right centres): one is the $H=1$ example whilst the other is $H=G\times G$. Equation \eqref{modinvsigma} becomes
$\cZ=b_1b_{G\times G}^t$ where $b_1$ resp.\ $b_{G\times G}$ are the matrices capturing the restrictions of the
previous two paragraphs.

For a type 1 example with a nontrivial $\psi$, consider $H=A_4\times A_4$, where $G=A_4$ is the alternating group of order 12 and $\omega=1$. Now, $H^2(A_4;\bbT)\cong \bbZ_2$, so let $\psi''$  be the nontrivial one.  $A_4$ has four conjugacy classes: those of (1), (12)(34), (123) and (132), with stabilisers $A_4,\bbZ_2\times\bbZ_2,\bbZ_3,\bbZ_3$ respectively. $\cD_1(A_4)$ has 14 sectors but most aren't relevant here. The map $\beta^{\psi''}$ is trivial on
the stabilisers $\bbZ_3$, but restricts to the unique nontrivial alternating bicharacter on $\bbZ_2^2$ (this follows from Theorem 2  in \cite{Hugh} or Lemma 1(b) in \cite{EG7}). We find that the sigma-restrictions of these two module categories are
slightly different:
\begin{align}\psi'=1:&\ \ [1,\mathbf{1}]\mapsto [1,\mathbf{1}]+[(12)(34),{}_{^{++}}]+[(123),\mathbf{1}]+[(132),\mathbf{1}]\,,\\
\psi'=\psi'':&\ \ [1,\mathbf{1}]\mapsto[1,\mathbf{1}]+[(12)(34),{}_{^{+-}}]+[(123),\mathbf{1}]+[(132),\mathbf{1}]\,,\end{align}
where we write $s_1s_2$ for the one-dimensional representation of $\bbZ_2^2$ sending
$(1,0)\mapsto s_1$ and $(0,1)\mapsto s_2$, and $\psi'$ is as in Corollary 2.

Consider next the lattice VOAs $\cV(A_2\oplus E_6)$, $\cV(A_8)$ and $\cV(L)$ where $L=\mathrm{span}\{D_7\oplus\sqrt{36}\bbZ,(\frac{1}{2},\frac{1}{2},\frac{1}{2},\frac{1}{2},\frac{1}{2};\frac{9}{\sqrt{36}})\}$. These are all
$\bbZ_3$ orbifolds of the holomorphic VOA $\cV(E_8)$, and their categories of modules correspond to the  three inequivalent twisted doubles of $\bbZ_3$, namely $\cD_{\omega_0}(\bbZ_3)$, $\cD_{\omega_1}(\bbZ_3)$ and $\cD_{\omega_2}(\bbZ_3)$ respectively using the notation of \eqref{cyclic3cocycle}. According to Theorem 1, there are
precisely two nontrivial extensions of $\cV(A_2\oplus E_6)$ (i.e. type 1 module categories of $\cD_{\omega_0}(\bbZ_3)$),
namely $[0,0]\oplus[1,1]\oplus[2,2]=\cV(E_8)$ ($H=1$, $\psi=1$) and  $[0,0]\oplus[1,2]\oplus[2,1]\cong\cV(E_8)$ ($H=\bbZ_3\times\bbZ_3$, $\psi=1$). Here we write 
$[0],[1],[2]$ for the simple modules of both $\cV(A_2)$ and $\cV( E_6)$, and denote by $[k,k']=[k]\otimes
[k']$ the simple modules of $\cV(A_2\oplus E_6)$.
According to Theorem 1, there is
precisely one nontrivial extension of $\cV(A_8)$ (i.e. type 1 module category of $\cD_{\omega_1}(\bbZ_3)$),
namely $[0]\oplus[3]\oplus[6]=\cV(E_8)$ ($H=1$, $\psi=1$). Here $[0],[1],\ldots,[8]$ denote the
simple modules of $\cV(A_8)$, enumerated so that fusion is addition mod 9. Similarly, there is only 1 for  $\cV(L)$, namely $H=1,\psi=1$, with $\oplus_{j=0}^2[0,\frac{9j}{\sqrt{36}}]\cong\cV(E_8)$.

Consider now the symmetric group $G=S_3$. Then $H^3(S_3;\bbT)\cong\bbZ_6$, with explicit 3-cocycles $\omega_k$ given in e.g. equation (6.20) of \cite{CGR}. In Table 1 we collect together the type 1 module categories for each of these 6 twisted doubles. In all cases, $\psi$ must be trivial. The column `MTC' is the MTC of the local modules, as always a twisted group double $\cD_{\omega_q}(F)$. Conveniently,  $H^2(K;\bbT)=0$ for all subgroups $K\le S_3$, which means all 2-cocycles $c^\omega_g$ are coboundary, so all characters are (projectively) equivalent to linear ones, and we can give all
sectors uniform names regardless of $\omega_k$.  In particular, enumerate the 8 sectors of $\cD_{\omega_k}(S_3)$ as $\chi_0=[1,\mathbf{1}]$, $\chi_1=[1,sgn],\chi_2=[1,\tau]$, $\chi_{3+k}=[(123),\xi^k]$, $\chi_{6+l}=[(12),(-1)^l]$, as in \cite{EP1}; label the $n^2$ sectors of $\cD_{\omega_q}(\bbZ_n)$
as $[i,\phi_j]$ where $0\le i,j<n$ and $\phi_j(1)=e^{2\pi ij/n}$.  In the column `branching rules', we list the sigma-restrictions of each
sector of $\cD_{\omega_q}(F)$ in order; for $F\cong \bbZ_n$ we order the sectors $[i,\phi_j]$  lexicographically.  We use the fact that,
for each prime $p$, the restriction of the $p$-primary part of $H^3(G;\bbT)$ to $H^3(P;\bbT)$ is injective, where
$P$ is a $p$-Sylow subgroup of $G$.
The module categories for $\omega=1$ were first given in \cite{EP1}.

\medskip\centerline {{\bf Table 1.} Type 1 module categories for $\cD_{\omega_k}(S_3)$}\medskip

\centerline{
{
\vbox{\tabskip=0pt\offinterlineskip
  \def\tablerule{\noalign{\hrule}}
  \halign to 5in{
    \strut#&
    \hfil#&\vrule#&\vrule#&\hfil#&\vrule#&    
\hfil#&\vrule#&\hfil#&\vrule#&\hfil#&\vrule#&\hfil#&\vrule#&\hfil#&\vrule#&\hfil#&\vrule#&\hfil#&\vrule#&\hfil#
\tabskip=0pt\cr
&$\omega_k\ \ \,$&\,&&\hfill\,\,$H$\,\,\hfill&&\hfill\,\,branching rules\,\,\hfill&&\hfill MTC\hfill\cr
\tablerule\tablerule&$\omega_0\ \ \ $&\,&&\hfill $1$\hfill&&\hfill $\,\,\chi_0+\chi_1+2\chi_2\,\,$\hfill&&\hfill Vec\hfill\cr
\tablerule&$\omega_0\ \ \ $&\,&&\hfill $\Delta(\bbZ_2)$\hfill&&\hfill $\,\chi_0+\chi_2,\,\chi_1+\chi_2,\,\chi_6,\,\chi_7\,$\hfill&&\hfill$\cD_{\omega_0}(\bbZ_2)$\hfill\cr
\tablerule&$\omega_0\ \ \ $&\,&&\hfill $\bbZ_2\times\bbZ_2$\hfill&&\hfill$\, \chi_0+\chi_2+\chi_6\,$ \hfill&&\hfill Vec\hfill\cr
\tablerule&$\omega_0$\ \ \ &\,&&\hfill $\,\,\Delta(\bbZ_3)\,$\hfill&&\hfill $\,\chi_0+\chi_1,\,\chi_2,\,\chi_2,\,\chi_3,\,\chi_4,\,\chi_5,\,\chi_3,\,\chi_4,\,\chi_5\,\,$\hfill&&\hfill $\cD_{\omega_0}(\bbZ_3)$ \hfill\cr
\tablerule&$\omega_0\ \ \ $&\,&&\hfill $\,\bbZ_3\times\bbZ_3\,$\hfill&&\hfill $\,\chi_0+\chi_1+2\chi_3\,$\hfill&&\hfill Vec\hfill\cr
\tablerule&$\omega_0\ \ \ $&\,&&\hfill $\,\Delta(S_3)\,$\hfill&&\hfill $\chi_i$\hfill&&\hfill$\cD_{\omega_0}(S_3)$\hfill\cr
\tablerule&$\omega_0$\ \ \ &\,&&\hfill $\,\Delta(S_3)(\bbZ_3\times 1)$\hfill&&\hfill $\chi_0+\chi_3,\,\chi_1+\chi_3,\,\chi_6,\,\chi_7$\hfill&&\hfill$\cD_{\omega_0}(\bbZ_2)$\hfill\cr
\tablerule&$\omega_0\ \ \ $&\,&&\hfill $S_3\times S_3$\hfill&&\hfill $\chi_0+\chi_3+\chi_6$\hfill&&\hfill Vec\hfill\cr
\tablerule\tablerule&$\omega_3\ \ \ $&\,&&\hfill 1\hfill&&\hfill $\,\,\chi_0+\chi_1+2\chi_2\,\,$\hfill&&\hfill Vec\hfill\cr
\tablerule&$\omega_3$\ \ \ &\,&&\hfill \,$\Delta(\bbZ_2)$\,\hfill&&\hfill $\,\chi_0+\chi_2,\,\chi_1+\chi_2,\,\chi_6,\,\chi_7\,$\hfill&&\hfill$\cD_{\omega_1}(\bbZ_2)$\hfill\cr
\tablerule&$\omega_3\ \ \ $&\,&&\hfill $\Delta(\bbZ_3)$\hfill&&\hfill $\,\chi_0+\chi_1,\,\chi_2,\,\chi_2,\,\chi_3,\,\chi_4,\,\chi_5,\,\chi_3,\,\chi_4,\,\chi_5\,$\hfill&&\hfill$\cD_{\omega_0}(\bbZ_3)$\hfill\cr
\tablerule&$\omega_3\ \ \ $&\,&&\hfill $\bbZ_3\times\bbZ_3$\hfill&&\hfill$\,\chi_0+\chi_1+2\chi_3\,$\hfill&&\hfill Vec\hfill\cr
\tablerule&$\omega_3$\ \ \ &\,&&\hfill $\Delta(S_3)$\hfill&&\hfill$\chi_i$\hfill&&\hfill$\cD_{\omega_3}(S_3)$\hfill\cr
\tablerule&$\omega_3\ \ \ $&\,&&\hfill $\Delta(S_3)(\bbZ_3\times 1)$\hfill&&\hfill$\chi_0+\chi_3,\,\chi_1+\chi_3,\,\chi_6,\,\chi_7$\hfill&&\hfill$\cD_{\omega_1}(\bbZ_2)$\hfill\cr
\tablerule\tablerule&$\omega_{\pm 2}\ \ \ $&\,&&\hfill 1\hfill&&\hfill$\,\,\chi_0+\chi_1+2\chi_2\,\,$\hfill&&\hfill Vec\hfill\cr
\tablerule&$\omega_{\pm 2}\ \ \ $&\,&&\hfill $\Delta(\bbZ_2)$\hfill&&\hfill$\,\chi_0+\chi_2,\,\chi_1+\chi_2,\,\chi_6,\,\chi_7\,$\hfill&&\hfill$\cD_{\omega_0}(\bbZ_2)$\hfill\cr
\tablerule&$\omega_{\pm 2}$\ \ \ &\,&&\hfill $\bbZ_2\times\bbZ_2$\hfill&&\hfill$\, \chi_0+\chi_2+\chi_6\,$\hfill&&\hfill Vec\hfill\cr
\tablerule&$\omega_{\pm 2}\ \ \ $&\,&&\hfill $\Delta(\bbZ_3)$\hfill&&\hfill$\,\chi_0+\chi_1,\,\chi_2,\,\chi_2,\,\chi_3,\,\chi_4,\,\chi_5,\,\chi_3,\,\chi_4,\,\chi_5\,$\hfill&&\hfill $\cD_{\omega_{\pm 1}}(\bbZ_3)$\hfill\cr
\tablerule&$\omega_{\pm 2}\ \ \ $&\,&&\hfill $\Delta(S_3)$\hfill&&\hfill$\chi_i$\hfill&&\hfill$\cD_{\omega_{\pm 2}}(S_3)$\hfill\cr
\tablerule\tablerule&$\omega_{\pm 1}\ \ \ $&\,&&\hfill 1\hfill&&\hfill$\,\,\chi_0+\chi_1+2\chi_2\,\,$\hfill&&\hfill Vec\hfill\cr
\tablerule&$\omega_{\pm 1}\ \ \ $&\,&&\hfill $\Delta(\bbZ_2)$\hfill&&\hfill$\,\chi_0+\chi_2,\,\chi_1+\chi_2,\,\chi_6,\,\chi_7\,$\hfill&&\hfill $\cD_{\omega_1}(\bbZ_2)$\hfill\cr
\tablerule&$\omega_{\pm 1}\ \ \ $&\,&&\hfill $\Delta(\bbZ_3)$\hfill&&\hfill$\,\chi_0+\chi_1,\,\chi_2,\,\chi_2,\,\chi_3,\,\chi_4,\,\chi_5,\,\chi_3,\,\chi_4,\,\chi_5\,$\hfill&&\hfill$\cD_{\omega_{\pm 1}}(\bbZ_3)$\hfill\cr
\tablerule&$\omega_{\pm 1}\ \ \ $&\,&&\hfill $\Delta(S_3)$\hfill&&\hfill$\chi_i$\hfill&&\hfill$\cD_{\omega_{\pm 1}}(S_3)$\hfill\cr
\noalign{\smallskip}}}}}

\section{{Proofs}}

By a \textit{transversal} $t_{\bar{k}}$ for $K/N$ we mean a choice of coset representatives $t_{\bar{k}}\in\overline{k}\subset K$. 
Write $\beta$ for the quantity $\beta^\psi$ defined in \eqref{betag}, $c_g$ for the
cocycle  $c_g^\omega$ of \eqref{2cocycle}, and define $\widetilde{\omega}$ by \eqref{tildeomega}.

There are two different sorts of extensions {in Theorem 1}, what we call \textit{type 1{$\Delta$}} (arising in the quantum Galois correspondence) and \textit{type 1{$N$}}. The
general extension is a combination of these. Type $1\Delta$ has $H=\Delta_K$ and $\psi=1$, for any subgroup $K\le G$ (no constraint comes from $\omega$). Type $1N$ has  $H=\Delta_G(1\times N)$ where $N\trianglelefteq G$ is normal, with $\psi(G,H)=1$, $\beta^\psi_{N\times 1}(1\times N)=1$ and $[\widetilde{\omega}]_H=[1]$.

\subsection{Type $1\Delta$}

We will study the general case through the type $1N$ module category $(\Delta_K(1\times N),\psi)$ of $\cD_\omega(K)$, transferred to $\cD_\omega(G)$ through the
type $1\Delta$ module category $(\Delta_K,1)$ of $\cD_\omega(G)$. For this reason, we need to first understand the type $1\Delta$ case.

Given any pair $(H,\psi)$, Proposition 1 associates the $\cD_\omega(G)$ module category  ${}_{\Delta_G}\cC(G^2)^{\tilde{\omega},\psi}_H$ and  full category  ${}_H^{\,\psi}\cC(G^2)^{\,\tilde{\omega},\psi}_H$. Any module category is associated
an associative algebra $A=\theta$ in $\cD_\omega(G)$; the module category is type 1 iff $A$ is commutative in the braided sense $\mu_A\circ c_{A,A}=\mu_A$.  
The general relation between (bi)modules and bundles that we need is worked out in
Lemma 3.1 of \cite{EP1}.

Consider first type $1\Delta$, so
$H=\Delta_K$. By $\cF(G/K)$ we mean the algebra of functions on the set $G/K$ of left cosets. The multiplication in this algebra is 
pointwise, $(f_1f_2)([g])=f_1([g])\,f_2([g])$ (hence is commutative). $\cF(G/K)$ carries a (true) representation of $G$, namely
Ind$_K^G\mathbf{1}$, obtained through the $G$-action on $G/K$. The desired algebra $A_\Delta=\theta_\Delta$, which must be a bundle in ${}_{\Delta_G}\cC(G^2)^{\,\tilde{\omega}}_{\Delta_G}$, is a copy of $\cF(G/K)$ on each point of $\Delta_G$,
with $\Delta_G^{L\times R}$ action given by \eqref{actionleft}-\eqref{actionmixed} and algebra multiplication $\mu_\Delta$ from equivariance \eqref{multequiv}, both applied to the fibre $\cF(G/K)$ at $(1,1)\in G^2$.
As a sector in $\cD_\omega(G)$, it is $[1,\mathrm{Ind}_K^G\mathbf{1}]$.
The indecomposable right $A_\Delta$-modules in ${}_{\Delta_G}\cC(G^2)_{\Delta_G}^{\,\tilde{\omega}}$ are precisely those bundles in ${}_{\Delta_G}\cC(G^2)^{\,\tilde{\omega}}_{\Delta_G}$ of the form $[(g,1), \mathrm{Ind}_{C_K(g)}^{C_G(g)}\overline{\chi}]$ where $\overline{\chi}\in\mathrm{Irr}_{c_g}(C_K(g))$; 
 this $A_\Delta$-module corresponds to the indecomposable  bundle $[(g,1),\overline{\chi}]\in{}_{\Delta_G}\cC(G^2)^{\,\tilde{\omega}}_{ \Delta_K}$. That correspondence defines the equivalence between the category of right $A_\Delta$-modules in ${}_{\Delta_G}\cC(G^2)_{\Delta_G}^{\,\tilde{\omega}}$ and the category ${}_{\Delta_G}\cC(G^2)^{\,\tilde{\omega}}_{\Delta_K}$.
 $A_\Delta$ is indeed commutative in the braided sense because the fibre $\cF(G/K)$ above the orbit representative $(1,1)$  is a commutative algebra in the classical sense. Thus  $(H,\psi)=(\Delta_K,1)$ is  type 1, as desired.

The above bijection also gives sigma-restriction, which is the functor from $A$-modules to objects in $\cD_\omega(G)$, forgetting the $A$-module structure. In full generality, it sends $[(g,g'),\overline{\chi}]\in{}_{\Delta_G}\cC(G^2)^{\,\tilde{\omega}}_{
\Delta_K}$, where $\overline{\chi}$ is a (projective) representation of the stabiliser ${}^{(g,g')}\Delta_K\cap\Delta_G\cong C_K(g^{\prime-1}g)$, to the bundle $[(g,g'),\mathrm{Ind}_{ {}^{(g,g')}\Delta_K\cap\Delta_G}^{{}^{(g,g')}\Delta_G\cap\Delta_G}\overline{\chi}]\in{}_{\Delta_G}\cC(G^2)^{\,\tilde{\omega}}_{
\Delta_G}$ (note that ${}^{(g,g')}\Delta_G\cap\Delta_G\cong C_G(g^{\prime-1}g)$). This coincides with the sigma-restriction
given in \cite{EG3} for the module category $(\Delta_K,1)$.

The local $A$-modules in ${}_{\Delta_G}\cC(G^2)_{\Delta_G}^{\,\tilde{\omega}}$ are precisely those whose sigma-restriction has a well defined
twist, or $T$-eigenvalue (see Theorem 3.2 of \cite{KO}). It is clear that any $A_\Delta$-module
of the form $[(k,1),\mathrm{Ind}_{C_K(k)}^{C_G(k)}\overline{\chi}]$ for $k\in K$ is local: $k\in Z(C_G(k))\cap Z(C_K(k))$ so any subrepresentation of $
\mathrm{Ind}_{C_K(k)}^{C_G(k)}\overline{\chi}$ has the same value of twist, namely $\overline{\chi}(k)/\overline{\chi}(1)$. In terms of the full system ${}_{\Delta_K}\cC(G^2)^{\,\tilde{\omega}}_{\Delta_K}$, this corresponds
to the bundle $[(k,1),\overline{\chi}]$. This class of local $A_\Delta$-modules are
precisely those which are also bundles over the subgroupoid ${}_{\Delta_K}\cC(K^2)^{\,\tilde{\omega}}_{\Delta_K}$; these form the MTC $\cD_\omega(K)$. The easiest way to see that these exhaust all local $A_\Delta$-bundles, is to compute dimensions: Theorem 4.5 of \cite{KO} tells us the global dimension of the MTC of local $A_\Delta$-bundles, which by definition is the square-root of the sum of
squares of the quantum-dimensions of all of its sectors, equals the global dimension of $\cD_\omega(G)$ (which is $|G|$) divided by the quantum-dimension of $A_\Delta$ in $\cD_\omega(G)$, which is $|G/K|$ (the quantum-dimension of a bundle in ${}_{\Delta_G}\cC(G^2)^{\,\tilde{\omega}}_{\Delta_G}$ is the dimension of its total space divided by $|G|$). Thus $\cD_\omega(K)$ has the same global dimension as the MTC of local $A_\Delta$-modules, of which it is a full subcategory, and hence the two must be equal. 
The map $\iota_!$ of Theorem 1(b)(iii) sends $[k,\overline{\chi}]$ to $[k,\mathrm{Ind}_{C_K(k)}^{C_G(k)}\overline{\chi}]$.

\subsection{Proof of type 1-ness in general case}

The nesting of type 1 module categories is developed in Proposition 4.16  of
\cite{FFRS} and Proposition 2.3.2 of \cite{Dav}; we summarise it here. Suppose $\cM_1$
is a type 1 module category for an MTC $\cC$, and let $\cC_1$ be its MTC of local modules.
Suppose $\cM_2$
is a type 1 module category for $\cC_1$, let $A_2$ be the associated commutative algebra, and let $\cC_2$ be its MTC of local modules. Then $A_2$ and indeed all of $\cC_2$ are objects in $\cC_1$ (this map is the forgetful functor called sigma-restriction). Likewise, every object in $\cC_1$ is an object in
$\cC$. Let $A$ be the algebra $A_2$ regarded as an object in $\cC$ in this sense. $A$  is an algebra
in $\cC$: the multiplication $\mu_A:A\otimes_{\cC} A\rightarrow A$ is the composition of the projection to the fusion product $A_2\otimes_{\cC_1} A_2$ in $\cC_1$, with the multiplication
$\mu_{A_2}:A_2\otimes_{\cC_1} A_2\rightarrow A_2$, which we identify with (sigma-restrict to)  $A$ ($V\otimes_{\cC_1} W$ is naturally a quotient of
$V\otimes_{\cC} W$ --- see Theorem 1.5 of \cite{KO}). The corresponding module category $\cM$ of $A$-modules in $\cC$ is the desired nesting or transfer of $\cM_2$ to $\cC$ via $\cM_1$. $A$ will be type 1 when both $\cM_1$ and $\cM_2$ are type 1. The MTC of local modules for $\cM$ is the sigma-restriction of $\cC_2$ first to $\cC_1$ then to $\cC$.

Consider now type $1N$, so  $H=\Delta_G(1\times N)$. The algebra $A_N=\theta_N$ here is the  bundle in ${}_{\Delta_G}\cC(G^2)^{\tilde{\omega},\psi}_H$ (hence  ${}_{\Delta_G}\cC(G^2)^{\,\tilde{\omega}}_{\Delta_G}$) consisting of  $\bbC$ attached to each point in $H$, with a $\Delta_G^L\times H^R$-action given by \eqref{actionleft}-\eqref{actionmixed}, where the fibre $\bbC_{(1,1)}$ has the trivial representation of the stabiliser $\Delta_G^{adj}$. The algebra multiplication $\mu_N:\bbC_{(1,1)}\otimes\bbC_{(n,1)}\rightarrow \bbC_{(n,1)}$ is that of $\bbC$; to obtain the product
$\mu_N:\bbC_{(gn,g)}\otimes\bbC_{(g'n',g')}\rightarrow\bbC_{(gg'n^{g'}n',gg')}$, use equivariance \eqref{multequiv} as before.
Then $A_N$ decomposes into a sum of $\cD_\omega(G)$ sectors as
$A_N=\sum_n [(n,1),\beta_{(n,1)}]$, where the sum is over representatives $n$ of each $G$-conjugacy class (i.e. $G^{adj}$-orbit)
in $N$. A  bundle in ${}_{\Delta_G}\cC(G^2)^{\,\tilde{\omega}}_{\Delta_G}$ is an $A_N$-module, iff it is also a bundle in ${}_{\Delta_G}\cC(G^2)^{\tilde{\omega},\psi}_{H}$. The indecomposable right $A_N$-modules in ${}_{\Delta_G}\cC(G^2)_{\Delta_G}^{\,\tilde{\omega}}$ are precisely the indecomposable bundles in ${}_{\Delta_G}\cC(G^2)^{\tilde{\omega},\psi}_{H}$ with $\Delta_G^L\times H^R$ equivariance restricted to $\Delta_G^{L\times R}$. Since the fibre of $A_N$ above $(1,1)\in G^2$ is a commutative algebra, this is also type 1.

Finally, consider the general case, $H=\Delta_K(1\times N)$ with $\psi$ as in Theorem 1. We know from last subsection that $(\Delta_K,1)$ is a type 1 (in fact type $1\Delta$) module category $\cM_1$ for 
$\cD_\omega(G)$, whose MTC of local modules is $\cD_\omega(K)$. We know from the previous paragraph that $(\Delta_K(1\times N),\psi)$ is a type 1 module category (in fact type $1N$) for $\cD_\omega(K)$. Nesting these, we obtain the module category $(\Delta_K(1\times N),\psi)$ of
$\cD_\omega(G)$, which is therefore type 1. The description of sigma-restriction last subsection
makes it easy to describe the algebra and MTC of local modules for this nesting. In particular, the
algebra $A_N$ as a bundle in ${}_{\Delta_K}\cC(K^2)_{\Delta_K}^{\,\tilde{\omega}}$ consists
of a copy of $\bbC$ above each point in $\Delta_K(1\times N)$; it is the identical bundle in ${}_{\Delta_K}\cC(G^2)_{\Delta_K}^{\,\tilde{\omega}}$; it corresponds to the bundle in ${}_{\Delta_G}\cC(G^2)_{\Delta_K}^{\,\tilde{\omega}}$ consisting of a copy of $\bbC$ (carrying the representation $\beta_{(1,n)}$) above every point in $\Delta_G(1\times N)$; and finally it corresponds to the bundle  in ${}_{\Delta_G}\cC(G^2)_{\Delta_G}^{\,\tilde{\omega}}$ described in section 3.1.

This shows that the pairs $(H,\psi)$ described in the second sentence of Theorem 1(a) are indeed type 1. The converse is Theorem 2.19 of \cite{Dav3}, though it also falls out from the 
treatment of general module categories for $\cD_\omega(G)$ which we develop in \cite{EG8}. The final sentence of Theorem 1(a) is deferred to after Claim 2, in section 4.7.

\subsection{{Proof of alpha-induction}}

In a type 1 theory, alpha-induction can be defined purely categorically, in terms of the
algebra $A$  \cite{KO}, as we sketched in section 2.1: we get a functor $\alpha^H$ from $\cD_\omega(G)$
to right $A_H$-modules ${}_{\Delta_G}\cC(G^2)^{\tilde{\omega},\psi}_H$, or functors 
$\alpha^H_\pm$ to $A_H$-$A_H$-bimodules ${}_{H}^{\,\psi}\cC(G^2)^{\tilde{\omega},\psi}_{H}$. 
The more fundamental is $\alpha^H$; we revisit $\alpha^H_\pm$ next subsection
when we compute the branching rules. Unlike sigma-restriction of local modules, alpha-induction for a nesting is not
simply a composition of the alpha-inductions for each component, but we have learned enough that it is now straightforward to do it in one step.

Consider any $(\Delta_K(1\times N),\psi)$ as in Theorem 1(a). As is clear from the previous two subsections, its commutative algebra is $A_H=\sum_n[(n,1),\mathrm{Ind}_{C_K(n)}^{C_G(n)}\beta^*_{(n,1)}]$, where the sum is over representatives of the $K^{adj}$ orbits in  $N$. Note that $A_\Delta=[(1,1),\mathrm{Ind}_K^G1]$ is a
subalgebra of $A_H$, and so an $A_H$-module is also an $A_\Delta$-module. In fact, a bundle
in ${}_{\Delta_G}\cC(G^2)^{\tilde{\omega},\psi}_H$ is an $A_H$-module iff it is an $A_\Delta$-module, hence expressible in the form $\sum [g,\mathrm{Ind}_{C_K(g)}^{C_G(g)}\overline{\chi}]$,
and also the corresponding bundle $\sum[(g,1),\overline{\chi}]\in{}_{\Delta_G}\cC(G^2)^{\,\tilde{\omega}}_{\Delta_K}$ is also a bundle in ${}_{\Delta_G}\cC(G^2)^{\tilde{\omega},\psi}_H$ (though with $\Delta_G^L\times H^R$-equivariance restricted to $\Delta_G^L\times \Delta_K^R$).
This association of any bundle in ${}_{\Delta_G}\cC(G^2)^{\tilde{\omega},\psi}_H$ to some bundle in ${}_{\Delta_G}\cC(G^2)^{\,\tilde{\omega}}_{\Delta_G}$ is sigma-restriction here.

Then alpha-induction sends a sector $[(g,1),\chi]$ to the bundle product $[(g,1),\chi]\otimes A_H$. Let's compute this from the intermediate point, i.e. as a bundle in ${}_{\Delta_G}\cC(G^2)^{\,\tilde{\omega}}_{\Delta_K}$\,: from \eqref{tensorofbundle} we obtain
\begin{align}\alpha^H([g,\chi])\leftrightarrow&\,\sum_n[(g,1),\chi]\otimes_{\Delta_G}[(n,1),\beta_{(n,1)}^*]\nonumber\\ =&\,\sum_n\sum_{h\in C_G(g)\backslash G/C_K(n)}[(g^hn,1), \mathrm{Ind}_{C_K(n)\cap C_K(g^hn)}^{C_K(g^hn)}(\chi^{hn}\,\beta^*_{(n,1)})]\,,\end{align}
where we use the $\Delta_G^L$ action to choose  more convenient orbit representatives. Interpreting this as a bundle in ${}_{\Delta_G}\cC(G^2)^{\tilde{\omega},\psi}_H$ is the same as dropping the sum over $n$, so we obtain \eqref{alphagen}.

As a consistency check, the algebra $A^H$ can be identified with the bundle
${\iota}=[(1,1),\mathbf{1}]$ in  ${}_{\Delta_G}\cC(G^2)_H^{\tilde{\omega},\psi}$, as discussed in \cite{Ev}. Its transpose $\overline{\iota}$ is the bundle $[(1,1),\mathbf{1}]$ in ${}_H^{\,\psi}\cC(G^2)_{\Delta_G}^{\,\tilde{\omega}}$. The sigma-restriction  functor from the full system ${}_H^{\,\psi}\cC(G^2)_H^{\,\tilde{\omega},\psi}$ to
$\cD_\omega(G)$, can be computed as the product   $\iota a\overline{\iota}$ of bundles.

\subsection{{Local modules and branching rules}}

In this subsection we prove Theorem 1(b)(ii), i.e. we work out explicitly sigma-restriction ${}_{\Delta_G}\cC(G^2)_H^{\tilde{\omega},\psi}\rightarrow {}_{\Delta_G}\cC(G^2)_{\Delta_G}^{\,\tilde{\omega}}$, and identify the local bundles in ${}_{\Delta_G}\cC(G^2)_H^{\tilde{\omega},\psi}$.

The type 1$\Delta$ case, i.e $(\Delta_K,1)$, is worked out in section 4.1. In particular,
the MTC of local modules is $\cD_\omega(K)$, and sigma-restriction Res$ :\mathrm{Ver}_\omega(K) \rightarrow \mathrm{Ver}_\omega (G)$  of local modules is  Res$([k , \overline{\chi}] )= [k, \mathrm{Ind}_{C_K({k})}^{C_G({k})} \overline{\chi} ]\,,$
 where $\overline{\chi}\in\mathrm{Irr}_{c_k}(C_K(k))$
and Ind is induction to a $c_k$-projective $G$-character.  Of course this coincides with the $K$-theoretic map $\iota_!$ where $\iota$ is the embedding $K\hookrightarrow G$, hence it coincides with Theorem 1(b)(iii) in this special case.

Now turn to the type $1N$ case, i.e. $H=\Delta_G(1\times N)$ and $\psi$ as in Theorem 1. 
Choose any indecomposable bundle $[(k,1),\chi]$ in ${}_{\Delta_G}\cC(G^2)_H^{\tilde{\omega},\psi}$. Its support is $\Delta_G(k,1)H=(\mathrm{cl}(k)\times 1)H$ where cl$(k)$ is the conjugacy class of $k$ in $G$. The stabiliser  of $(k,1)$ is St$_H(k,1)=\{(h_+,h_-)\in H:h_-k=kh_+\}\cong C_H(k)$. As we know, this bundle sigma-restricts to a bundle in ${}_{\Delta_G}\cC(G^2)_{\Delta_G}^{\,\tilde{\omega}}$ simply by restricting $ \Delta_G^L\times H^R$ equivariance
to $\Delta_G^{L\times R}$. If we write the set $\mathrm{cl}(k)N$ as a disjoint union $\cup_n\mathrm{cl}(kn)$ over certain $n\in N$ (possible, since $N$ is normal in $G$), this sigma-restriction becomes $\sum_n[(kn,1),\chi']$, where $\chi'$ is the appropriate character (identified shortly) of the stabiliser St$_H(kn,1)\cap\Delta_G\cong C_G(kn)$. 

To identify that character $\chi'$, let $V$ be the fibre above $(k,1)$ in the bundle $[(k,1),\chi]$, carrying a $C_H(k)$-irrep $\rho$ realising the character $\chi$. Choose any $n\in N$ and let $V'$ be the fibre above
$(kn,1)$ in that bundle; we can identify $V'$ with $V$ through the invertible map $v\mapsto v'=v.(n,1)$. Note that  $(g,g')\in \mathrm{St}_H(kn,1)$ iff
$({}^ng,g')\in \mathrm{St}_H(k,1)$. We can compare 
$\rho'(g,g')v':=(g',g').(v'.(g,g^{\prime})^{-1})\in V'$ to $\rho({}^ng,g')v:=(g',g').(v.({}^ng,g')^{-1})$ through the moves $(g',g').\bigl((v.(n,1)).(g,g^{\prime})^{-1}\bigr)\approx
(g',g').\bigl(v.(ng^{-1},g^{\prime\,{-1}})\bigr)=(g',g').\bigl(v.(({}^ng,g')^{-1}(n,1))\bigr)\approx
\bigl((g',g').(v.({}^ng,g')^{-1})\bigr).(n,1)$, where `$\approx$' means equality up to a phase from \eqref{actionleft}-\eqref{actionmixed}; keeping track of the phases,  we obtain 
\begin{equation}\rho'(g,g')v'=\beta_{(n,1)}({}^ng^{-1},g^{\prime-1})^*
\omega(k,n,g^{-1})^*\omega(k,{}^ng^{-1},n)\,\omega(g',g^{\prime-1}k,n)\bigl(\rho({}^ng,g')v\bigr)'\label{rho'}\end{equation}
Again, $\rho'$ is the representation of St$_H(kn,1)$ in the fibre above $(kn,1)$ in the bundle $[(k,1),\chi]$; we are interested in its restriction $\chi'$ to St$_H(kn,1)\cap\Delta_G\cong C_G(kn)$.

We want to identify the local $A^H$-modules, i.e. the indecomposable bundles $[(k,1),\chi]$ in ${}_{\Delta_G}\cC(G^2)_H^{\tilde{\omega},\psi}$ which, when we restrict to a bundle in ${}_{\Delta_G}\cC(G^2)_{\Delta_G}^{\,\tilde{\omega}}$, has a well-defined twist.  Choose any $n\in N$ and  $g\in C_G(kn)$; then $(g,g)\in \mathrm{St}_H(kn,1)$. Writing $\Delta_h=(h,h)$, 
the twist (or $T$-matrix eigenvalue) of the indecomposable bundle $[(k',1),\chi']$ in ${}_{\Delta_G}\cC(G^2)^{\,\tilde{\omega}}_{\Delta_G}\cong \cD_\omega(G)$ is
$\chi'(k')/\chi'(1)$. So the bundle $[(k,1),\chi]$ in  ${}_{\Delta_G}\cC(G^2)_H^{\tilde{\omega},\psi}$ is local iff $\chi(k,k)/\chi(1,1)$ has modulus 1 (i.e. $\rho(k,k)$ is scalar), and 
\begin{equation}\label{local}\chi(k,k)=\beta_{(n,1)}(k^{-1}n^{-1},n^{-1}k^{-1})^*\omega(k,n,n^{-1}k^{-1})^*\omega(k,k^{-1}n^{-1},n)\,\omega(kn,n^{-1},n)\,\chi(nk,kn)\end{equation} for all $n\in N$. Note that some characters $\chi$ of St$_H(k,1)$ may not satisfy $|\chi(k,k)|=\chi(1,1)$, since $(k,k)$ will not in general lie in the centre of
St$_H(k,1)$.

Sigma-restriction and local modules for general type 1 $(H,\psi)$ is now a straightforward
nesting of the treatments for types $1\Delta$ and $1N$, and is as in Theorem 1(b)(ii).

\subsection{The associated modular tensor category}

In this subsection we identify the MTC of local modules as a twisted group double.

The type 1$\Delta$ case was worked out in section 4.1, so consider first type ${1N}$, i.e. $H=\Delta_G(1\times N)$ and  $\psi$ as in Theorem 1. We will find that the effect of the 2-cocycle $\psi$ is surprisingly subtle.  
Write $\overline{G}:= G/N$ and $\pi:G\rightarrow\overline{G}$ for the corresponding homomorphism.
 In the following, we  identify
$C_{\Delta_G}(ng,1)$ with $C_G(ng)$ and  $C_H(g)/N$ with $C_{\overline{G}}(gN)$.

We begin with the elementary observation that a projective representation of a group $K$ defines a projective representation for $G/N$, $N$ the
projective kernel:

\medskip\noindent\textbf{Claim 1.} \textit{Let $N$ be a normal subgroup 
of $K$, and fix any transversal $t_{\bar{k}}$ for $K/N$. Let $c\in Z^2(K;\bbT)$ be normalised.}

\smallskip\noindent\textbf{(a)} \textit{Suppose $\chi\in\mathrm{Irr}_c(K)$ satisfies $|\chi(n)|=\chi(1)$ for all $n\in N$. Then $\phi(n):=\chi(n)/\chi(1)$ satisfies}
\begin{equation}\phi(n)\phi(n')=c(n,n')\phi(nn')\,,\ \phi(n^k)\,c(n,k)=\phi(n)\,c(k,n^k)\ \forall n,n'\in N\,,\ k\in K\,.\label{claim1conds}\end{equation} 

\noindent\textbf{(b)} \textit{Suppose $\phi:N\rightarrow\bbT$ satisfies \eqref{claim1conds}, and define 
\begin{equation}\nonumber\overline{c}(\overline{h},\overline{k}):=\phi(t_{\bar{h}}t_{\bar{k}}t_{\bar{h}\bar{k}}^{-1})\,c(t_{\bar{h}},t_{\bar{k}})\,c(t_{\bar{h}}t_{\bar{k}}t_{\bar{h}\bar{k}}^{-1},t_{\bar{h}\bar{k}})^*\,.\end{equation}
 Then
the relation  $\overline{\chi}(\overline{k})=\chi(t_{\bar{k}})$ defines 
 a  bijection between the set of $\chi\in\mathrm{Irr}_{c}(K)$ satisfying $\chi(n)=\phi(n)\chi(1)$ for all $n\in N$, and $\overline{\chi}\in\mathrm{Irr}_{\bar{c}}(K/N)$.}\medskip 

\noindent\textit{Proof.} Part (a) is trivial: let $\rho$ be the projective representation associated to $\chi$, and compute $\rho(n)\rho(n')$ and $\rho(k)\rho(n)$.

To prove (b), first let $\chi$ be as above, again realised by  $\rho$. First note that 
 the matrices $\rho(n)$   are scalar for $n\in N$  (this follows from
the triangle inequality applied to $\chi(n)=\phi(n)\chi(1)$). Now define $\overline{\rho}(\overline{h})
=\rho({t_{\overline{h}}})$. Then for all $\overline{h},\overline{k}\in {K/N}$,
$$\overline{\rho}(\overline{h})\overline{\rho}(\overline{k})=c(t_{\bar{h}},t_{\bar{k}})\rho(t_{\bar{h}}t_{\bar{k}}) 
={c}(t_{\bar{h}},t_{\bar{k}})\,c(t_{\bar{h}}t_{\bar{k}}t^{-1}_{\bar{h}\bar{k}},t_{\bar{h}\bar{k}})^*\rho(t_{\bar{h}}t_{\bar{k}}t^{-1}_{\bar{h}\bar{k}})\rho(t_{\bar{h}\bar{k}})
=\bar{c}(\bar{h},\bar{k})\bar{\rho}(\bar{h}\bar{k})$$
as desired. Conversely, given $\overline{\chi}\in
 \mathrm{Irr}_{\bar{c}}(K/N)$, define $\chi(nt_{\bar{k}})=c(n,t_{\bar{k}}) \phi(n)\overline{\chi}(\overline{k})$ for all $n\in N,\overline{k}\in K/N$,
 then the same calculation shows $\chi\in\mathrm{Irr}_c(K)$. Furthermore, these two maps $\chi\mapsto\overline{\chi}$
 and $\overline{\chi}\mapsto\chi$ are inverses.  \qquad\textit{QED to Claim 1}\medskip

Consider any sector $[g,\chi]\in\mathrm{Ver}_\omega(G)$ appearing in the sigma-restriction of a local bundle in ${}_{\Delta_G}\cC(G^2)_H^{\tilde{\omega},\psi}$. Then Claim 1 with $K$ there replaced with $C_H(g)$, $c$ replaced with $c_g$, 
and $\phi$ obtained from \eqref{local}, tells us that the $\chi\in\mathrm{Irr}_{c_g}(C_H(g))$
satisfying \eqref{local} are in natural bijection with the  $\overline{\chi}\in\mathrm{Irr}_{\bar{c}_g}(C_{\overline{G}}(\overline{g}))$, where $\bar{c}_g=\overline{c}_g^\omega\,\overline{c}_g^\psi$ for
\begin{align}\label{omegapart}\overline{c}_g^\omega(\overline{h},\overline{k})=&\,\omega(ng, n, g)^*\omega(n, ng, g) \,\omega(n, ng^2, g^{-1}n^{-1})\, c_g(n, g)\,c_g(t_{\bar{h}},t_{\bar{k}})\,c_g(n,t_{\bar{h}\bar{k}})^*\\ \label{psipart}
\overline{c}_g^\psi(\overline{h},\overline{k})=&\,\beta_{(g,g)}(n,1)\end{align}
and where we write $n=t_{\bar{h}}t_{\bar{k}}t_{\bar{h}\bar{k}}^{-1}$.

Note that the 2-cocycle $\psi$ and the 3-cocycle $\omega$ are completely independent:
the pair $(H,\psi)$ is type 1 for $\cD_\omega(G)$, iff $(H,1)$ is type 1 for $\cD_\omega(G)$; moreover, in that case it is also type 1 for $\cD_1(G)$. In addition, in \eqref{omegapart},\eqref{psipart} we factorised the multiplier $\overline{c}_g=\overline{c}_g^\omega\overline{c}_g^\psi$ where
$\overline{c}_g^\omega$ depends on $\omega$ but not $\psi$, and $\overline{c}_g^\psi$
depends on $\psi$ but not $\omega$.

We want to find a 3-cocycle  on $\overline{G}$ which is responsible (in the
sense of \eqref{2cocycle}) 
for these multipliers $\overline{c}_g$. The preceding paragraph means that it (if it exists) will factorise into a part $\overline{\omega}$ seeing only $\omega$, and a part
$\overline{\omega}^\psi$ seeing only $\psi$. 

Fix a transversal $t_{\bar{g}}$ for $G/N$. Recall that  $\widetilde{\omega}|_H$ must be coboundary.

\medskip\noindent\textbf{Proposition 2.} \cite{Dav3,Jon2} \textit{Let  $F$ be a  2-cochain on $H$ such that
$dF=\widetilde{\omega}|_H$.
 Define
\begin{align}\nonumber\overline{\omega}(\overline{g},\overline{h},\overline{k})=
\omega\bigl(k_1,t_{\bar{h}},t_{\bar{k}}\bigr)\,&
F((n_1,1),(n_{2},1))\,F((n_3,1),(n_2n_1,1))\times\\ \label{baromega} &F((n_3^{-1},1),(n_3,1)) \,\beta^F_{(k_1,k_1)}(n_4,1)\\
\overline{\omega}^\psi(\overline{g},\overline{h},\overline{k})=\psi((n_2,1),(n_1,1))&\,{\psi((n_4^{-1},1),(n_5,1))^*}\,\beta_{(t_{\bar{g}},t_{\bar{g}})}(n_4^{-1},1)\,,
\label{omegaJones} 
\end{align}
where ${k}_1=t_{\bar{g}\bar{h}\bar{k}}t^{-1}_{\bar{k}}t^{-1}_{\bar{h}}$, $n_{1}=t_{\bar{g}\bar{h}}t_{\bar{k}}t_{\bar{g}\bar{h}\bar{k}}^{-1}$, $n_{2}=t_{\bar{g}}t_{\bar{h}}t_{\bar{g}\bar{h}}^{-1}$, $n_4=t_{\bar{h}}t_{\bar{k}}t_{\bar{h}\bar{k}}^{-1}$, $n_{3}={}^{{k}_{1}}(n_4^{-1})$,  $n_5=t_{\bar{g}} t_{\bar{h}\bar{k}} t_{\bar{g}\bar{h}\bar{k}}^{-1}$, and where $\beta^F$ is as in \eqref{betag}. Then both $\overline{\omega},
\overline{\omega}^\psi\in Z^3(\overline{G};\bbT)$ and are normalised, and when substituted into \eqref{2cocycle} recover $\overline{c}_g^\omega$ and $\overline{c}_g^\psi$ respectively.}
\medskip

Choosing a different transversal $t_{\star}$ won't change the cohomology class of the 3-cocycles.
Up to notational differences, \eqref{baromega} agrees with equation (20) in \cite{Dav3} when the 2-cocycle $\psi$ is taken to be 1. The calculation of $\overline{\omega}^\psi$ is also contained in \cite{Dav3} (when $\omega=1$), but since we will later need to exploit more structure (see especially section 4.9)  we prefer the treatment in \cite{Jon2}, even though it necessitates chopping the calculation of the 3-cocycle into the two parts. 
The treatments are equivalent, indeed the explicit calculations given in the Theorem 2.17 proof would make equal sense in \cite{Jon2}. It would be nice though to understand conceptually why the two contexts are equivalent in this sense. The argument that these 3-cocycles recover the appropriate 2-cocycles
follows from \cite{Dav3}.

 The source of \eqref{omegaJones}
 is the thesis of Vaughan Jones \cite{Jon2}. Jones verified in  section 4.2 of \cite{Jon2},
 that for any  2-cocycle $\mu\in Z^2(N;\bbT)$ and map
 $\lambda:G\times N\rightarrow \bbT$ satisfying the following equations
 \begin{align}
 \lambda(m,n)=&\,\mu(m,n^m)\,{\mu(n,m)^*}\,,\label{jones1}\\
 \lambda(gh,n)=&\,\lambda(g,n)\,\lambda(h,n^g)\,,\label{jones2}\\
 \lambda(g,mn)\,{\lambda(g,m)^*}\,{\lambda(g,n)^*}=&\,\mu(m,n)\,{\mu(m^g,n^g)^*}\,,\label{jones3}\\
 \lambda(1,n)=\lambda(g,1)=&\,\mu(1,n)=\mu(n,1)=1\,,\label{jones4}\end{align}
for all $m,n\in N,g,h\in G$, the map $\overline{\omega}^\psi$ 
 is indeed a 3-cocycle on $\overline{G}$. Now, $\lambda(g,n)=\beta'_g(n)$ and $\mu(m,n)=\psi'(m,n)$ satisfy these equations: \eqref{jones1}-\eqref{jones3} is
\eqref{betag}-\eqref{bichar2} respectively, and \eqref{jones4} follows because $\psi$ is normalised, so
we can use Jones' result.
 
Jones was interested in  determining complete invariants for actions of  finite groups $G$ on the
 hyperfinite II$_1$ factor. His invariants consist of a normal subgroup $N$, and $\lambda,\mu$ satisfying these equations. A nice review of this circle of ideas is \cite{RSW}. Incidentally, in section 4.9 we need the converse, that any solution $\lambda,\mu$ comes from a $\psi$ in this way. So this point of contact between our context and that of Jones is pronounced.

Inflating, in the usual sense, an
indecomposable bundle in ${}_{\Delta_{\bar{G}}}\cC(\overline{G}^2)_{\Delta_{\bar{G}}}^{\,\tilde{\omega}'}$ (where we write $\omega'=\overline{\omega}\overline{\omega}^\psi$)  gives a (local) indecomposable bundle in ${}_{\Delta_G}\cC(G^2)_{\Delta_G}^{\,\tilde{\omega}}$. The equivalence of the category of local modules in ${}_{\Delta_G}\cC(G^2)^{\tilde{\omega},\psi}_{H}$ and $\cD_{\bar{\omega}'}(\overline{G})$ for the appropriate $\overline{\omega}'$ is established in Theorem 2.17 of \cite{Dav3}.

Sigma-restriction Ver$_{\bar{\omega}\bar{\omega}^\psi}(\overline{G})\rightarrow \mathrm{Ver}_\omega (G)$ is the map $\pi^*$ of Theorem 1(b)(iii)
corresponding to the group homomorphism $\pi:G\rightarrow \overline{G}$. We can read off from the proof of Claim 1(b) that 
an indecomposable bundle $ [\overline{g},\overline{\chi}]\in\mathrm{Ver}_{\bar{\omega}\bar{\omega}^\psi}(\overline{G})$ is sent to
\begin{equation}[\overline{g},\overline{\chi}]\mapsto \sum_{{g}}\bigl[{g},{\widetilde{\phi}}\,\overline{\chi}\circ\pi\bigr]\,,
\label{sigma2i}\end{equation}
where  the sum over $g$ is over representatives of the distinct conjugacy classes in $G$ projecting to the $\overline{G}$-conjugacy 
class of $\overline{g}$, and $\widetilde{\phi}(h)=c_g(ht^{-1}_{hN},t_{hN})\,\phi(ht^{-1}_{hN})$ where $\phi$ is
as in \eqref{local}. As explained in the previous paragraph, the quantities 
${\widetilde{\phi}}\,\overline{\chi}\circ\pi_{{g}}$
can be identified with $c_g^\omega$-projective characters of $C_G(g)$ behaving  appropriately with respect to $n\in C_N(g)$, as
 \eqref{local} requires.

The general type 1 module category $(\Delta_K(1\times N),\psi)$ is a combination of the
type 1$\Delta$ pair $(\Delta_K,1)$ followed by the type 1$N$ pair  $(\Delta_K(1\times N),\psi)$. Then the MTC of 
local bundles for the desired type 1 module category, as well as sigma-restriction from
local bundles to Vec$_\omega(G)$, is obtained from the type 1$N$ and 1$\Delta$ ones as explained
in the beginning of section 4.2. The result is given in Theorem 1(b)(iii).

\subsection{Proof of Corollary 1}

Take $G=\langle g\rangle\cong\bbZ_n$. Write $m=|K|$ and $m'=|N|$. Then $m'|m|n$ and $K=\langle
g^{n/m}\rangle$, $N=\langle g^{n/m'}\rangle$, and $H=\Delta_K(1\times N)\cong \bbZ_{m'}\times
\bbZ_{m}$, where the $\bbZ_{m'}$ factor has generator $g_1=(g^{n/m'},0)$ and the $\bbZ_m$ factor
has generator $g_2=(g^{n/m},g^{n/m})$.

Explicit $k$-cocycles are worked out in \cite{HWY} for all abelian groups and all $k$. 
For the group $H\cong\bbZ_{m'}\times\bbZ_m$, $H^3(H;\bbT)\cong \bbZ_{m'}^2\times\bbZ_m$
with representatives
\begin{align}\omega_{a_1,a_2,a_{12}}(g_1^{i_1}g_2^{i_2},g_1^{j_1}g_2^{j_2},g_1^{k_1}g_2^{k_2})=&
\exp\left(2\pi i \frac{a_1i_1}{m'}\left[\frac{j_1+k_1}{m'}\right]\right)\times\nonumber\\ &\exp\left(2\pi i \frac{a_2i_2}{m}\left[\frac{j_2+k_2}{m}\right]\right)
\exp\left(2\pi i \frac{a_{12}i_2}{m'}\left[\frac{j_1+k_1}{m'}\right]\right)\nonumber\end{align}
for integers $0\le a_1,a_{12}<m'$ and $0\le a_2<m$.
Now $\widetilde{\omega}_q|_H$ is the 3-cocycle given by
$$\widetilde{\omega}_q(g_1^{i_1}g_2^{i_2},g_1^{j_1}g_2^{j_2},g_1^{k_1}g_2^{k_2})=\exp\left(2\pi i \frac{q}{m}(\frac{m}{m'}i_1+i_2)\left[\frac{j_1+k_1}{m'}+\frac{j_2+k_2}{m}\right]\right)\exp\left(-2\pi i \frac{q}{m}i_2\left[\frac{j_2+k_2}{m}\right]\right)$$
This is required to be coboundary, so in particular the restriction to the subgroup $\langle g_1\rangle \cong \bbZ_{m'}$ must
be coboundary, but that restriction is clearly $q\omega_1'$ where $\omega_1'$ is the generator of $H^3(\bbZ_{m'};\bbT)\cong\bbZ_{m'}$ given in \eqref{cyclic3cocycle}. This restriction is coboundary iff $m'|q$.
Write $[\tilde{\omega}_1]|_H=[\omega_{a_1,a_2,a_{12}}]$. Since the restriction of $\tilde{\omega}_1|_H$ to $\Delta_K$ is coboundary (in fact identically 1), $a_2=0$. Therefore $[\tilde{\omega}_1]|_H=[\omega_{a_1,0,a_{12}}]$ has order dividing $m'$, so  $[\tilde{\omega}_q]|_H=q[\tilde{\omega}_1]|_H
=[1]$. Hence the only conditions on $H$ (i.e. on $m',m$) are $m'|m|n$ and $m'|q$.

By the K\"unneth formula, $H^2(H;\bbT)\cong \bbZ_{m'}$. We find (e.g. \cite{HWY}) that $\psi$ can be taken to be
$$\psi_a(g_1^{i_1}g_2^{i_2},g_1^{j_1}g_2^{j_2})=\exp\left(ai_1j_2/m'\right)$$ 
where $0\le a<m'$.  Since $H^2(\bbZ_m;\bbT)=0$, the restriction of any $\psi_a$ to $\Delta_K$ is automatically coboundary. For the remaining condition on $\psi$, we compute
\begin{align}\psi_a((\frac{n}{m'}i_1,0),(0,\frac{n}{m'}j_2))=\,&\psi_a(g_1^{i_1},g_1^{-j_1}g_2^{mj_2/m'})=\exp\left(2\pi i\,ami_1j_2/(m')^2\right)\nonumber\\
\psi_a((0,\frac{n}{m'}j_2),(\frac{n}{m'}i_1,0))=\,&\psi_a(g_1^{-j_1}g_2^{mj_2/m'},g_1^{i_1})=1\nonumber\end{align}
These must be equal for all $i_1,j_2$, so we obtain $(m')^2|am$.

\subsection{{Proof of Corollary 2 and other things}}

\medskip\noindent\textbf{Claim 2.} \textit{Let $H=\Delta_K(1\times N)$.}

\medskip\noindent\textbf{(a)} \textit{Each class in  $Z^2(H;\bbT)$ contains a cocycle $\psi$ satisfying}
\begin{equation} \psi((gn,g),(km,k))=\beta^{\psi}_{(k,k)}(n,1)^*\,\psi((n^k,1),(m,1))\,\psi((g,g),(k,k))\,.\label{psiguaged}\end{equation}

\noindent{\textbf{(b)} \textit{Conversely, suppose $\psi'\in Z^2({1\times N};
\bbT)$, $\psi\in Z^2({\Delta_K};\bbT)$, and  $\beta_{(g,g)}(n,1)\in \bbT$ is a function 
on $g\in K,n\in N$, satisfying \eqref{bichar}: i.e. $\beta_{(gk,gk)}(n,1)=\beta_{(g,g)}(n,1)\,
\beta_{(k,k)}(n^g,1)$ for all $g,k\in K$, $n\in N$. Extend $\psi'$ to $H\times H$ by \eqref{psiguaged}; then $\psi'\in Z^2(H;
\bbT)$.}}\medskip

\noindent\textit{Proof.} Consider first part (a). Let $\psi'\in Z^2(H;\bbT)$ be arbitrary. Define all $f(k,k)=1$,
choose values $f(n,1)$ arbitrarily, and define $f(kn,k)=\psi'((k,k),(n,1))^*\,f(n,1)$. Then $\psi:=\psi'\,\delta f^*$
is a cocycle cohomologous to $\psi'$. We compute $\psi((g,g),(m,1))=1$, so \eqref{betag} collapses to $\psi((n,1),(k,k))=\beta^\psi_{(k,k)}(n,1)^*$ and the cocycle condition \eqref{twcoc} tells us 
$\psi((g,g),(km,k))=\psi((g,g),(k,k))$. 
Using this and \eqref{twcoc},\eqref{bichar2} we get 
$$\psi((n,1),(km,k))=\psi((n,1),({}^{k}m,1))\,\psi((n\,{}^km,1),(k,k))\,\psi(({}^km,1),(k,k))^*$$
$$=\frac{\psi((n,1),({}^km,1))}{\beta^\psi_{(k,k)}({}^km,1)^*}\left(\frac{\beta^\psi_{(k,k)}(n,1)\beta^\psi_{(k,k)}({}^km,1)\psi((n,1),({}^km,1))}{\psi((n^k,1),(m,1))}\right)^*$$
which collapses to agree with \eqref{psiguaged}. Using this and the cocycle condition, one finds \eqref{psiguaged} holds in general. Part (b) is also straightforward.
\qquad\textit{QED to Claim 2}\medskip

Suppose $(H,\psi)$ is type 1 (so $H=\Delta_K(1\times N)$), and has trivial left and right centres $\cC_l(A)=\cC_r(A)=\mathrm{Vec}_\bbC$.  Theorem 4.5 of \cite{KO} tells us that the global
dimension of the MTC of local modules equals the global dimension of $\cD_\omega(G)$, namely $|G|$, divided by the quantum-dimension of $A$, namely $|N|\,|G/K|$, and so equals $|K|/|N|$. But Vec$_\bbC$ has
dimension 1, which forces $N=K$.  
Assume $\psi$ is put in the form \eqref{psiguaged}, which we can rewrite as
\begin{equation} \psi((gn,g),(km,k))=\psi((n,1),(k,k))\,\psi'(n^k,m)\,\psi((g,g),(k,k))\,,\end{equation}
for some $\psi'\in Z^2(K;\bbT)$.
 The condition $\psi(K,H)=1$ is equivalent to the requirement that $\psi((g,g),(k,k))=1$ for all $g,k\in K$. The requirement on $\psi((n,1),(1,n'))$ is equivalent to the condition that $\psi((n,1),(k,k))=\psi'(k{}^{-1},n)\,
 \psi'(n^{k},k{}^{-1})^*$. We recover the form for $\psi$ given in Corollary 2. We also require that the class $[\widetilde{\omega}|_H]$ be trivial. But the cocycle $\omega(k_1,k_2,k_3)\,\omega(k_1',k_2',k_3')^*$
 is coboundary on $G^2$ iff $\omega$ is coboundary on $G$. This proves one direction of Corollary 2.
 
 For the other direction, let $\psi'\in Z^2(K;\bbT)$. Define $\psi(K,K)=1$ and $\psi((n,1),(k,k))$ by the formula given in the previous paragraph. Then according to Claim 2(b), $\psi$ defined as
 in Corollary 2 is a 2-cocycle of $H=K^2$, provided $\beta_{(gk,gk)}(n,1)=\beta_{(g,g)}(n,1)\,\beta_{(k,k)}(n^g,1)$, i.e. provided
 \begin{equation}\psi'(k^{-1}g^{-1},n)\,\psi'(n^g,g^{-1})\,\psi'(n^{gk},k^{-1})=\psi'(n^{gk},k^{-1}g^{-1})\,
 \psi'(g^{-1},n)\,\psi'(k^{-1},n^g)\,.\nonumber\end{equation}
 But this follows from the relations $\psi'(n^{gk},k^{-1})=\psi'(n^gk,k^{-1})\,\psi'(k^{-1},n^g)\,\psi'(k^{-1},n^gk)^*$, $\psi'(n^{gk},k^{-1}g^{-1})=\psi'(n^gk,k^{-1}g^{-1})\,\psi'(k^{-1},g^{-1}n)\,\psi'(k^{-1},n^gk)^*$, $\psi'(n^{g}k,k^{-1}g^{-1})=\psi'(n^gk,k^{-1})\,\psi'(n^g,g^{-1})\,\psi'(k^{-1},g^{-1})^*$, and
 $\psi'(k^{-1},g^{-1}n)=\psi'(k^{-1}g^{-1},n)\,\psi'(k^{-1},g^{-1})\,\psi'(g^{-1},n)^*$, coming from the cocycle condition \eqref{2coc}.
 
 We can now conclude the proof of the final  sentence of Theorem 1(a): choose $\psi$ to be the
  cocycle in its cohomology class whose existence is promised by Claim 2(a). We see directly from
 \eqref{psiguaged} that $\psi(\Delta_G,H)=1$. The condition $\beta^\psi_{N\times 1}(1\times N)=1$,
 which is equivalent to $\psi((n,1),(1,n'))=\psi((1,n'),(n,1))$, holds for $\psi$ iff it holds for all cocycles
 in $[\psi]$.
  
\subsection{Proof of Theorems 2 and 4(a)}

Since $G\le S_k$, $(\cA^{\otimes k})^G$ contains $(\cA^{\otimes k})^{S_k}$, and the 3-cocycle $\omega_G$ for $(\cA^{\otimes k})^G$ is the restriction of that for 
$(\cA^{\otimes k})^{S_k}$. So it suffices to prove Theorem 2 for the full symmetric group $G=S_k$.

First, let us recall some generalities about group cohomology (see e.g. \cite{Brown}). 
Given an embedding $\iota_K:K\hookrightarrow G$, we get restriction
$\iota_K^*:H^n(G;\bbT)\rightarrow H^n(K;\bbT)$, a group homomorphism.
Given any finite group $G$ and prime $p$ dividing its order, it is elementary that restriction $\iota_P^*$ to a $p$-Sylow subgroup $P$ sees only the $p$-primary part of $H^n(G;\bbT)$ (i.e. the elements annihilated by some power of $p$), i.e. it kills the $p'$-primary part for any prime $p'\ne p$. Less obvious is that \textit{transfer} (see e.g. Theorem 10.3 in \cite{Brown}) says that restriction $\iota_P^*$ of the $p$-primary part of $H^n(G;\bbT)$ to $P$ is injective. 

Nakaoka  (Theorem 5.8 of \cite{Nak}) has proved that for any embedding $\iota:S_m\hookrightarrow S_k$ for $m\le k$, the map
$\iota^*$ is surjective, with kernel which is a direct summand.  Moreover (his Corollary 6.7), for $n<k/2$
restriction $\iota^*:H^n(S_k;\bbT)\rightarrow H^n({S_{k-1}};\bbT)$ is an isomorphism. The first few groups $H^3(S_k;\bbT)$ are
$$\bbZ_2\,,\ \bbZ_6\,,\ \bbZ_{12}\times \bbZ_2\,,\ \bbZ_{12}\times \bbZ_2\,,\ \bbZ_{12}\times \bbZ_2\times\bbZ_2\,,
$$
for $S_2$, $S_3$, $S_4$, $S_5$, resp. $S_k$ for $k\ge 6$, where it stabilises.

In particular, for $k\ge 3$ the 3-torsion subgroup in $H^3(S_k;\bbT)$ is $\bbZ_3$; combining Nakaoka's result with transfer allows us to identify (via $\iota^*$) that 3-torsion in $H^3(S_k;\bbT)$ ($k\ge 3$) with $H^3(C_3;\bbT)$ where
$C_3=\langle(123)\rangle\cong\bbZ_3$. We write $[\omega^{(3)}_q]\in H^3(S_k;\bbT)$ for the unique 3-torsion class which restricts to $[\omega_q]\in H^3(C_3;\bbT)$ in \eqref{cyclic3cocycle}. This is the cocycle appearing in the statement of Theorem 2.

\medskip\noindent\textbf{Claim 3.} \textit{Suppose some class $[\omega]\in H^3(S_k;\bbT)$ is nontrivial. Then there is either a cyclic subgroup 
$C\le S_k$ of order $\le 4$ such that the restriction of $[\omega]$ to $C$ is nontrivial, or $k\ge 6$ and the restriction of
$\omega$ to  $\langle (12)\rangle\times\langle(34)\rangle\times\langle(56)\rangle\cong \bbZ_2\times\bbZ_2\times\bbZ_2$ is cohomologous to $\omega_{iii}$ of \eqref{omegaiii}.}\medskip

\noindent\textit{Proof.} Without loss of generality we may assume $k=6$, by stability. Write (noncanonically) $H^3(S_6;\bbT)=\langle x\rangle\times\langle y\rangle\times\langle z\rangle$ where $x,y,z$ have orders $12,2,2$ respectively.

 Suppose first that $[\omega]\in
H^3(S_k;\bbT)$ has order a multiple of 3. Then by the above discussion  it must restrict to an order 3 element of $H^3(C_3;\bbT)$, where $C_3$ is as above.

The 2-primary part of $H^3(S_k;\bbT)$ is more delicate.  Suppose first that $[\omega]\in
H^3(S_k;\bbT)$ has order  4 (so $[\omega]$ is of the form $x^{\pm 3}y^iz^j$). Then by Nakaoka it must restrict to an order-4 class in  $H^3(S_4;\bbT)$, hence by transfer to an order-4 class for its Sylow 2-subgroup $\langle (1324),(12)\rangle\cong D_4$, the dihedral group of order 8. Explicit cocycles can be computed e.g. by the homological algebra library HAP in GAP; a cocycle of order 4 in $H^3(D_4;\bbT)\cong\bbZ_4\times\bbZ_2
\times \bbZ_2$ is
$$\omega_{(4)}(a^lb^m,a^{l'}b^{m'},a^{l''}b^{m''})=\left\{\begin{matrix}e^{\pi i(-1)^ml/2}&\mathrm{if}\ m'=0\ \mathrm{and}\ 2l'+2l''+m'+m''\ge 8\\ e^{\pi i(-1)^{m+1}l/2}&\mathrm{if}\ m'=1\ \mathrm{and}\ 2l''+m''>2l'+m'\\ 1&\mathrm{otherwise}\end{matrix}\right.$$
where we write $D_4=\langle a,b\,|\,a^4=b^2=(ab)^2=1\rangle$, so $0\le l<4$, $0\le m<2$, etc. We see  that the restriction of $\omega_{(4)}$ to $\langle a\rangle\cong \bbZ_4$ (i.e. choosing $m=m'=m''=0$) is $\omega_1$ (recall \eqref{cyclic3cocycle}), hence order 4. The order-4 classes for $D_4$ equal $[\omega_{(4)}^{\pm 1}]$  times an order-2 class; in all cases it restricts to
an order-4 cocycle on  $C_4=
\langle (1234)\rangle\cong \bbZ_4$.

Of course, if $[\omega]=[\omega_{(4)}]^2$, then it is also nontrivial (in fact order 2) on $C_4$. So the only classes remaining to consider are $[\omega]\in\{y,yz,z\}$.

It is convenient to change scalars $\bbT$. The short exact sequence $0\rightarrow\bbZ\rightarrow\bbR\rightarrow\bbT\rightarrow 1$ gives a long exact sequence
in cohomology whose connecting homomorphisms $\alpha:H^n(G;\bbT)\rightarrow H^{n+1}(G;\bbZ)$ are isomorphisms.
The short exact sequence $0\rightarrow \bbZ\buildrel \times 2\over \rightarrow \bbZ\buildrel\pi\over\rightarrow\bbZ_2\rightarrow 0$, where $\pi$ is reduction mod 2, gives the long exact sequence
\begin{equation}\label{longbeta}\cdots \rightarrow H^n(G;\bbZ)\buildrel \times 2\over \rightarrow H^n(G;\bbZ)\buildrel\pi\over\rightarrow H^n(G;\bbZ_2)\buildrel\beta\over\rightarrow H^{n+1}(G;\bbZ)\rightarrow\cdots\end{equation}
The sequence $0\rightarrow\bbZ_2\buildrel \epsilon\over \rightarrow\bbT\buildrel square\over \rightarrow\bbT\rightarrow 1$ where $\epsilon$ identifies $\bbZ_2$ with $\pm 1$, gives a long exact sequence 
\begin{equation}\label{longgamma}\cdots \rightarrow H^n(G;\bbZ_2)\buildrel \epsilon\over \rightarrow H^n(G;\bbT)\buildrel square\over\rightarrow H^n(G;\bbT)\buildrel\gamma\over\rightarrow H^{n+1}(G;\bbZ_2)\rightarrow\cdots\end{equation}
Finally, $0\rightarrow \bbZ_2\rightarrow\bbZ_4\buildrel\times 2\over \rightarrow\bbZ_2\rightarrow 0$ gives
a long exact sequence where the connecting homomorphism is called the Steenrod square $Sq^1$. All of these connecting maps are essentially the same: more precisely, $Sq^1=\pi\circ\beta=\gamma\circ\epsilon$. However,  $Sq^1$ is the most accessible: it is a derivation on the (commutative) ring $H^*(G;\bbZ_2)$ (where the product is the cup product), and $Sq^1(a)=a^2$ for any class of degree 1. In all cases, naturality of the long exact sequence implies that the restrictions
of cohomology groups intertwine the connecting maps.

A fairly complete description of the $\bbZ_2$-cohomology of the symmetric groups is provided in Chapter VI of \cite{AM}. As rings, we have the isomorphisms $H^*(S_2;\bbZ_2)\cong\bbZ_2[\sigma'_1]$, 
and $H^*(S_6;\bbZ_2)\cong \bbZ_2[\sigma_1,\sigma_2,\sigma_3,c_3]/(c_3(\sigma_3+\sigma_1\sigma_2))$, where the subscript indicates the degree of the generators. Restriction from
$S_6$ to $S_2\times S_2\times S_2=\langle (12),(34),(56)\rangle$ (where $S_6$ is the permutation group for $\{1,2,3,4,5,6\}$) sends $c_3$ to 0, and $\sigma_1,\sigma_2,\sigma_3$ to the generators of the same name of
a polynomial subalgebra $\bbZ_2[\sigma_1,\sigma_2,\sigma_3]$ of $H^3(S_2^3;\bbZ_2)\cong H^3(S_2;\bbZ_2)\otimes_{\bbZ_2} H^3(S_2;\bbZ_2)\otimes_{\bbZ_2} H^3(S_2;\bbZ_2)$. Moreover,   $Sq^1(\sigma_2)=
\sigma_1\sigma_2+\sigma_3+c_3$, $Sq^1(\sigma_3)=(c_3+\sigma_3)\sigma_1$, $Sq^1(c_3)=0$, and as always $Sq^1(\sigma_1)=\sigma_1^2$ and $Sq^1(ab)=Sq^1(a)\,b+a\,Sq^1(b)$. 
 
Equation \eqref{longbeta} tells us that $\beta(H^3(S_6;\bbZ_2))$ equals the order 2 classes in $H^4(S_6;\bbZ)$, namely Span$\{ 6\alpha(x),\alpha(y),\alpha(z)\}\cong\bbZ_2^3$, and that $\pi$ kills  $6\alpha(x)$.
Therefore we can identify (through $\beta$) the image of $Sq^1=\pi\circ\beta$ on $H^3(S_6;\bbZ_2)$ with the space Span$\{ \alpha(y),\alpha(z)\}\cong \bbZ_2^2$. We compute this image to be Span$\{Sq^1(\sigma_1^3),Sq^1(\sigma_3)\}$. This means we can identify $y=\epsilon(\sigma_1^3)$ and
$z=\epsilon(\sigma_3)$.

Restriction to $S_2^3$ of $\sigma_1^3$ and $\sigma_3$ (hence of $y,z$) is now straightforward,
and we see that $z$ restricts to 0  on $\langle (12)\rangle=:C_2$, while $y$ (hence $yz$) restricts to  $\epsilon(\sigma_1^{\prime\, 3})\ne 0$  on $C_2$. Using the interpretation of the image of restriction in terms of symmetric invariants, as discussed in section III.4 of \cite{AM} --- $S_3$ acts by permuting the three $S_2$'s --- we identify the restriction of $z$ with the 3-cocycle
$\omega_{iii}$, as desired. \textit{QED to Claim 3}\medskip

So it suffices to consider permutation orbifolds of the form  $(\cA^{\otimes k})^C$, where $C$ is a cyclic group generated by an order $n$ cycle $\pi$ in $S_k$, where $n=2,3,4$, as well as $(\cA^{\otimes 6})^{S_2\times S_2\times S_2}$.

We can read off the conformal weights of any simple module in $(\cA^{\otimes n})^{\langle\pi\rangle}$ from the work \cite{BHS,Ba,BDM,LX,KLX} on permutation orbifolds (or the \eqref{permorbchar} paragraph). For example, {\cite{LX}, Theorem 6.3e} gives the conformal weights of the simple summands of the restriction to that permutation orbifold of the unique $\pi$-twisted module of the conformal net $\cA^{\otimes n}$: these are $\frac{j}{n}+\frac{n^2-1}{24n}c$ where $c\in 8\bbZ_{\ge 0}$
is the central charge of $\cA$ and $0\le j<n$. (We thank Marcel Bischoff for correspondence on this point.)  Hence the conformal weights of the permutation orbifold with
$G\cong\bbZ_n$ for $n=2,4$  lie in $\frac{1}{n}\bbZ$. As was discussed at the end of section 2.2, this means that the twist $\omega$ occurring in those $n=2,4$ orbifolds is trivial. However, when $n=3$, 
the conformal weights  lie in $\frac{1}{9}\bbZ$; in particular $c/9$ (mod 1) is a conformal weight.
Again from the calculations at the end of section 2.2, we see this requires the specific twist $\omega=\omega_{c\,\, (\mathrm{mod}\,\, 3)}\in H^3(\bbZ_3;\bbT)$.

Putting this together with transfer and the results of Nakaoka, we identify the twist
$[\omega]$ of the permutation orbifold $(\cA^{\otimes k})^G$ as the restriction to  $H^3(G;\bbT)$
of the class we call $[\omega^{(3)}_{c\,\, (\mathrm{mod}\,\, 3)}]$. 

Let's turn to $(\cA^{\otimes 6})^{S_2\times S_2\times S_2}$, where $S_2\times S_2\times S_2=\langle(12)\rangle\times\langle(34)\rangle\times\langle(56)\rangle$. This clearly
equals  $((\cA\otimes\cA)^{S_2})^{\otimes 3}$. But we've learned the 3-cocycle for $(\cA\otimes\cA)^{S_2}$ will be trivial, so so will that for the $S_2\times S_2\times S_2$ orbifold.

This concludes the proof of Theorem 2.\medskip
 
 The first step in proving Theorem 4(a) is to show that Conjecture 2 holds here. Because $G$ is solvable, $(\cV^{\otimes k})^G$ is completely rational, and Mod$((\cV^{\otimes k})^G)$ is a modular tensor category $\cC$. This (more precisely, \textit{any} orbifold by $G$ of a holomorphic VOA) is the situation studied in sections 4 and 5 of \cite{Kir}, and we begin by reviewing what is obtained there. Let $A_G$ be the commutative algebra in $\cC$ capturing the extension $\cV^G\subset\cV$. Then  the isomorphism classes of simple objects in Mod$_\cC(A_G)$ are in natural bijection with $g\in G$. Choose a representative  $X_g$ for each $g$. Then for each $g,h\in G$ there are isomorphisms $\mu_{g,h}\in\mathrm{Hom}(X_g\otimes_{A_G}X_h,X_{gh})\in\bbC$ ($\otimes_{A_G}$ denotes the tensor product in the fusion category Mod$_\cC(A_G)$), and these will be unique up to nonzero constants. Computing morphisms
in Hom$(X_g\otimes_{A_G} X_h\otimes_{A_G} X_k,X_{ghk})\cong \bbC$ by introducing brackets in two ways, we can define $\omega_G(g,h,k)\in\bbC^\times$ by $\mu_{g,hk}\circ(1\otimes_{A_G}\mu_{h,k})=\omega_G(g,h,k)\,
\mu_{gh,k}\circ (\mu_{g,h}\otimes_{A_G} 1)$. Then $\omega_G$ defines a unique class in $H^3(G,\bbC^\times)\cong H^3(G,\bbT)$.

If $H\leq G$, then any $A_G$-module $X_h$ for $h\in H$, i.e. any $h$-twisted $\cV^{\otimes n}$-module, is trivially also an $A_H$-module. So the 3-cocycle $\omega_H$ is the restriction
of $\omega_G$. This means, as in the proof for Theorem 2, we can assume $G=S_k$ for $k\ge 6$, so we are in the case considered by Claim 3. The possibility that $[\omega]$ is detected by a cyclic
group of order $\le 4$  violates for example the character calculations given around \eqref{permorbchar}, or the conformal weight calculation  implicit in equation (3.6) of  \cite{BDM}. As before, the cocycle
on $(\cA^{\otimes 6})^{S_2\times S_2\times S_2}=((\cA\otimes\cA)^{S_2})^{\otimes 3}$ must also be trivial. We conclude that $[\omega]=1$.

By the main theorem of \cite{Kir}, we obtain $\cC\cong\cD_1(G)$, concluding the proof of Theorem 4(a). Note that the only place we used solvability of $G$ is in learning $(\cV^{\otimes k})^G$ is completely rational, so we obtain more generally that Mod$((\cV^{\otimes k})^G)\cong \cD_1(G)$ whenever $24|c$ and we know $(\cV^{\otimes k})^G$ is completely rational.

\subsection{Proof of Theorems 3 and 4(b)}

Fix any  finite group $G$ and normal subgroup $N$. Jones (see section 1.2 of \cite{Jon2}) defines a group $\Lambda(G,N)$, his \textit{characteristic invariant}, as follows. If the pairs $(\lambda,\mu)$ and $(\lambda',\mu')$ both satisfy \eqref{jones1}-\eqref{jones4},
then so does their product $(\lambda\lambda',\mu\mu')$. Let $Z$ denote the resulting (abelian) group 
of solutions $\lambda:G\times N\rightarrow\bbT$, $\mu\in Z^2(G;\bbT)$ to those equations. 
For any function $\eta:N\rightarrow\bbT$ with $\eta(1)=1$, and define $\lambda_\eta(g,n)=\eta(n)\,\eta(n^g)^*$ and $\mu_\eta(n,n')=\eta(nn')\,\eta(n)^*\eta(n')^*$. Then $\lambda_\eta,\mu_\eta$ are readily
seen to satisfy  \eqref{jones1}-\eqref{jones4}, and so the set of all of them form a subgroup of $Z$. Jones' group
$\Lambda(G,N)$ is defined to be the quotient of $Z$ by that subgroup.

Likewise, let $Z'$ be the set of all $\psi\in Z^2(H;\bbT)$ satisfying both $[\psi]|_{\Delta_G}=[1]$ and $\beta_{N\times 1}(1\times N)=1$ (the conditions on a type 1 $\psi$ in Theorem 1(a). Since $\beta^{\psi\psi'}=\beta^\psi\beta^{\psi'}$, $Z'$ is a subgroup of $Z^2(H;\bbT)$. Since $(n,1)$ and $(1,n')$ commute, $\beta^\psi_{(n,1)}(1,n')=1$ when 
$\psi$ is coboundary, so $Z'$ contains all 2-coboundaries. Define $\cG(G,N)$ to be the quotient of $Z'$ by all 2-coboundaries; it is the subgroup of $H^2(H;\bbT)$ consisting of all classes $[\psi]$ where
$(H,\psi)$ is type 1.

\medskip\noindent\textbf{Claim 4.} \textit{The map sending $(\lambda,\mu)$ to $\psi((gn,g),(kn',k)):=\mu(n^k,n')/\lambda(k,n)$ defines an isomorphism between the groups $\Lambda=\Lambda(G,N)$ and $\cG=\cG(G,N)$, with inverse sending $\psi$ to $\lambda(g,n)=\beta^\psi_{(g,g)}(n,1)$ and $\mu(n,n')=\psi((n,1),(n',1))$.}\medskip

\noindent\textit{Proof.} Call a 2-cocycle $\psi\in Z^2(H;\bbT)$ \textit{well-chosen} if it satisfies $\psi(\Delta_G,H)=1$ and $\beta^\psi_{N\times 1}(1\times N)=1$. Let $(H,\psi)$ be type 1. Then Theorem 1(a) says $\psi$ is cohomologous to a well-chosen one. Suppose 
$\psi$ and $\psi'$ are both well-chosen and cohomologous, and write $\psi'=\psi \,\delta f$ for some $f:H\rightarrow\bbT$ with $f(1,1)=1$. Because $\psi|_{\Delta_G}=\psi'|_{\Delta_G}$, we know
$f(gk,gk)=f(g,g)\,f(k,k)$. Because $\psi((g,g),(n,1))=\psi'((g,g),(n,1))$, we know $f(gn,g)=f(g,g)\,f(n,1)$.
We compute $(\delta f)((gn,g),(km,k))=f(n^km,1)\,f(n,1)^*f(m,1)^*$. Therefore $\psi,\psi'$ are both well-chosen and cohomologous iff there is a function $\eta:N\rightarrow\bbT$ with $\eta(1)=1$ and
$\psi'((gn,g),(km,k))\psi((gn,g),(km,k))^*=\eta(n^km)\,\eta(n)^*\eta(m)^*=:\psi_\eta((gn,g),(km,k))$. We've shown $\cG$ is the group $Z''$ of all well-chosen $\psi$, quotient the subgroup of all $\psi_\eta$.

Because $\mu\in Z^2(N;\bbT)$ and $\lambda$ satisfies \eqref{jones2}, Claim 2(b) tells us $\psi((gn,g),(km,k)):=\mu(n^k,m)/\lambda(k,n)$ lies in $Z^2(H;\bbT)$. The resulting $\psi$ is well-chosen: $\psi(\Delta_G,H)=1$ by \eqref{jones4}, and $\beta^\psi_{(n,1)}(1,n')=1$ by \eqref{jones1} and \eqref{jones2}. Moreover, $(\lambda_\eta,\mu_\eta)\mapsto \psi_\eta$, so passing to the quotients, we get a group homomorphism $F:\Lambda\rightarrow\cG$.
 
For $\psi$ well-chosen, define $\lambda(g,n)=\beta^\psi_{(g,g)}(n,1)=\psi((n,1),(g,g))^*$ and $\mu(n,n')=\psi((n,1),(n',1))$. As mentioned in section 4.5, $(\lambda,\mu)$ lie in $Z$.
 This map is the inverse of the map $Z\rightarrow Z''$ of the previous paragraph.
 The map sends $\psi_\eta$ to $(\lambda_\eta,\mu_\eta)$, so passing to the quotients we get the inverse $F^{-1}:\cG\rightarrow\Lambda$.  \textit{QED to Claim 4}\medskip

By an extension of $G$ by $N$, we mean a group $\widetilde{G}$ with a normal subgroup $\widetilde{N}\cong N$ for which $\widetilde{G}/\widetilde{N}\cong G$. Associated to the projection $p:\widetilde{G}\rightarrow G$ is \textit{inflation} $p^*:H^*(G;\bbT)\rightarrow H^*(\widetilde{G};\bbT)$. Associated to the embedding $\iota:\widetilde{N}\hookrightarrow 
\widetilde{G}$ is \textit{restriction} $\iota^*:H^*(\widetilde{G};\bbT)\rightarrow H^*(\widetilde{N};\bbT)$.

Jones  (amongst others) showed (Proposition 4.2.5 in \cite{Jon2}) that the standard 5-term restriction-inflation exact sequence for finite group extensions
$1\rightarrow \widetilde{N}\rightarrow \widetilde{G}\rightarrow {G}\rightarrow 1$ can be extended to the right:
\begin{align}
1\rightarrow &\,H^1(G;\bbT)\,{\buildrel p^*\over\rightarrow} \,H^1(\widetilde{G};\bbT)\,{\buildrel \iota^*\over\rightarrow} \,H^1(\widetilde{N};\bbT)^{{G}}\rightarrow H^2({G};\bbT)\,{\buildrel p^*\over\rightarrow} \,H^2(\widetilde{G};\bbT)\nonumber\\ &\,\,{\buildrel \iota^*\over\rightarrow} \,\Lambda(G,N)
\;{\buildrel\gamma\over\rightarrow}\;H^3({G};\bbT)\;{\buildrel p^*\over\rightarrow}\; H^3(\widetilde{G};\bbT)\label{JonesSequence}\end{align}
In fact this sequence can be continued indefinitely \cite{HJPR}, but we only need those next three terms.
The map $\gamma$ here coincides with the assignment \eqref{omegaJones}; by Claim 4 we can replace $\Lambda$ with $\cG$. The map $p^*$ there is inflation. 

It is well-known that a projective representation of a finite group $G$ can be lifted 
to a true (i.e. linear) representation of a finite extension of $G$. The cohomological fact underlying this
is that there is an extension $\widetilde{G}$ of $G$ (by the Schur multiplier) such that inflation
$H^2(G;\bbT)\rightarrow H^2(\widetilde{G};\bbT)$ is the 0 map. We need the $H^3$ analogue, which
is Lemma 7.1.2 in \cite{Jon2}. In fact this holds more generally:

\medskip\noindent\textbf{Proposition 3.} \cite{Op} \textit{Let $G$ be finite, and $q\ge 2$.
Then there exists an extension $\widetilde{G}$ by a finite abelian group (depending on $q$)  such that inflation
$H^q(G;\bbT)\rightarrow H^q(\widetilde{G};\bbT)$ is the 0 map.}\medskip

Finally, turn to the proof of Theorem 3. 
Choose any finite group $G$ and any 3-cocycle $\omega
\in Z^3(G;\bbT)$, and define $\widetilde{\omega}\in Z^3(G^2;\bbT)$ by \eqref{tildeomega}.
Choose an extension $\widetilde{G}$ of $G$ by some group $N$ so that inflation $p^*:H^3(G;\bbT)
\rightarrow H^3(\widetilde{G};\bbT)$ is trivial. Then exactness of \eqref{JonesSequence} says $\gamma$ must be surjective. This means we can choose a type 1 $\psi\in Z'$ so that $\overline{\omega}^\psi\in
H^3(G;\bbT)$ is cohomologous with $\omega$.

Fix any holomorphic conformal net $\cA$ in central charge 24 --- e.g. the Monstrous moonshine module. Use Cayley's Theorem to fix some embedding of $\widetilde{G}$ in some symmetric group $S_k$. Next, perform the permutation orbifold   $(\cA^{\otimes k})^{\tilde{G}}$. By Theorem 2,
this has category of representations $\cD_1(\widetilde{G})$. Then for $N$ and $\psi$ chosen as in the
previous paragraph, we know $(\Delta_{\tilde{G}}(1\times N),\psi)$ is type 1. Let $\cA'$ denote the
corresponding extension. Then by Theorem 1, $\cA'$ will be a completely rational conformal net with category of representations Rep$(\cA')$ equivalent to $\cD_\omega(G)$.

Of course, the extension $\widetilde{G}$ in Proposition 3 is solvable iff $G$ is, since $N$ is abelian here.
Since Conjecture 1 is known to hold for any solvable $G$, we get unconditionally that any $\cD_\omega(G)$ with $G$ solvable, is realised as the MTC Mod$(\cV)$ of some completely rational VOA
$\cV$.

After we completed this proof and shared it with Victor Ostrik, he later communicated an alternate proof of Theorem 3, which assumes our Theorem 2, which we now sketch. First, upgrade the homomorphism $K\rightarrow G$ to a tensor functor
$F: \mathrm{Vec}({K})\rightarrow \mathrm{Vec}_\omega(G)$. 
Let $I$ be the right adjoint of $F$ and let $A=I(1)$. Then $A$ has a natural lift
to the untwisted double of $K$, where it is an \'etale algebra.
Moreover the category Rep$_{\mathrm{Vec}({K})}(A)$ is tensor equivalent to Vec$_\omega(G)$.
The category of local $A$-modules is  
the centre of Rep$_{\mathrm{Vec}({K})}(A)$, giving the desired result.

The proof of Theorem 4(b) is identical, thanks to Theorem 4(a).

\newcommand\biba[7]   {\bibitem{#1} {#2:} {\sl #3.} {\rm #4} {\bf #5,}
                    {#6 } {#7}}
                    \newcommand\bibx[4]   {\bibitem{#1} {#2:} {\sl #3} {\rm #4}}

\def\ASENS            {Ann. Sci. \'Ec. Norm. Sup.}
\def\AM   {Acta Math.}
   \def\AnM              {Ann. Math.}
   \def\CMP              {Commun.\ Math.\ Phys.}
   \def\IJM              {Internat.\ J. Math.}
   \def\JAMS             {J. Amer. Math. Soc.}
\def\JFA              {J.\ Funct.\ Anal.}
\def\JMP              {J.\ Math.\ Phys.}
\def\JRA              {J. Reine Angew. Math.}
\def\JSP              {J.\ Stat.\ Physics}
\def\LMP              {Lett.\ Math.\ Phys.}
\def\RMP              {Rev.\ Math.\ Phys.}
\def\RNM              {Res.\ Notes\ Math.}
\def\RIMS             {Publ.\ RIMS.\ Kyoto Univ.}
\def\Inv              {Invent.\ Math.}
\def\npbp             {Nucl.\ Phys.\ {\bf B} (Proc.\ Suppl.)}
\def\nupb             {Nucl.\ Phys.\ {\bf B}}
\def\nup              {Nucl.\ Phys. }
\def\nupp             {Nucl.\ Phys.\ (Proc.\ Suppl.) }
\def\adma             {Adv.\ Math.}
\def\coma             {Con\-temp.\ Math.}
\def\PAMS             {Proc. Amer. Math. Soc.}
\def\PJM              {Pacific J. Math.}
\def\ijmp             {Int.\ J.\ Mod.\ Phys.\ {\bf A}}
\def\jpa              {J.\ Phys.\ {\bf A}}
\def\PLB              {Phys.\ Lett.\ {\bf B}}
\def\RIMS             {Publ.\ RIMS, Kyoto Univ.}
\def\Top               {Topology}
\def\TAMS             {Trans.\ Amer.\ Math.\ Soc.}

\def\Duke              {Duke Math.\ J.}
\def\K                 {K-theory}
\def\JOP               {J.\ Oper.\ Theory}

\vspace{0.2cm}\addtolength{\baselineskip}{-2pt}
\begin{footnotesize}
\noindent{\it Acknowledgement.}

The authors thank the University of Alberta Mathematics Dept, Cardiff School of Mathematics,  University of Warwick 
Mathematics Institute, Swansea University Dept of Computer Science and the  Isaac Newton Institute for Mathematical Sciences, Cambridge during the programme Operator Algebras: Subfactors and their applications, for generous hospitality while researching this
paper. They also benefitted greatly from Research-in-Pairs held at Oberwolfach and BIRS,
and the von Neumann algebra trimester at the Hausdorff Institute. 
DEE thanks  WWU M\" unster and Wilhelm Winter whilst TG thanks Karlstads Universitet and J\"urgen Fuchs for stimulating environments while part of this paper was written. We also thank Marcel Bischoff, Michael M\"uger and Victor Ostrik for comments.
Their research was supported in part by  EPSRC grant nos EP/K032208/1 and EP/N022432/1, PIMS, NSERC and SFB 878.

\end{footnotesize}


\begin{thebibliography}{99}
    \begin{scriptsize}

    \addcontentsline{toc}{section}{References}

    \setlength{\parskip}{-1ex}

\bibx{AM} {Adem, A., Milgram, R. J.} {Cohomology of Finite Groups.} (second edition) Berlin: Springer-Verlag,  2004.



\bibx{Bnt} {Bantay, P.} {Orbifolds, Hopf algebras, and the moonshine.} Lett. Math. Phys. 22 (1991), 187--194.

\bibx{Ba} {Bantay, P.} {Characters and modular properties of permutation orbifolds.} Phys. Lett. B419 (1998), 175--178.


\bibx{BDM} {Barron, K., Dong, C., Mason, G.} {Twisted sectors for tensor product vertex operator algebras associated to permutation groups.} {Commun. Math. Phys.} 227 (2002), 349--384.

\bibx{Bis} {Bischoff, M.} {Conformal net realizability of Tambara--Yamagami categories and generalized metaplectic modular categories.} arXiv://1803.04949.


\bibx{BE1}
{B\"ockenhauer, J., Evans, D. E.}
{Modular invariants, graphs and alpha-induction for
nets of subfactors. I.} {\CMP \,{197}, {361--386} {(1998)}; II. \CMP \, {200},
{57--103} {(1999)}; III. \CMP \, {205}, {183--228} {(1999)}}. 

\bibx{BE4}
{B\"ockenhauer, J., Evans, D. E.}
{Modular invariants from subfactors: Type I coupling
matrices and intermediate subfactors.} \CMP\,{213}, {267--289} {(2000)}.


\bibx{BEK1}
{B\"ockenhauer, J., Evans, D. E., Kawahigashi, Y.}
{On $\alpha$-induction, chiral generators and modular
invariants for subfactors.} \CMP\,{208}, {429--487} {(1999)}.

\bibx{BEK2}
{B\"ockenhauer, J., Evans, D. E., Kawahigashi, Y.}
{Chiral structure of modular invariants for subfactors.}
\CMP\,{210}, {733--784} {(2000)}.

\bibx{BHS} {Borisov, L., Halpern, M. B., Schweigert, C.} {Systematic approach to cyclic orbifolds.} {Int. J. Mod. Phys.} A13 (1998), 125--168.


\bibx{Brown} {Brown, K, S.}
{Cohomology of groups.} 
{Graduate Texts in Mathematics, 87. New York: Springer-Verlag,  1994.}


\bibx{CM} {Carnahan, S., Miyamoto, M.} {Regularity of fixed-point vertex operator subalgebras.} (arXiv:1603.05645v3 [math.RT]), 2016.


\bibx{CKLW} {Carpi, S., Kawahigashi, Y.,  Longo, R.,  Weiner, M.} {From vertex operator algebras to conformal nets and back.} Memoirs Amer. Math. Soc. (to appear); arXiv:1503.01260. 



\bibx{Che} {Cheng, C.} {A character theory for projective representations of finite groups.} {Linear Alg. Appl.} {469} (2015), 230--242.


\bibx{CGR} {Coste, A., Gannon, T., Ruelle, P.} {Finite group modular data.}
Nucl. Phys. B581 (2000), 679--717.


\bibx{CKM} {Creutzig, T., Kanade, S.,  McRae, R.} {Tensor categories for vertex operator superalgebra extensions.} arXiv:1705.05017.


\bibx{Dav} {Davydov. A.} {Modular invariants for group-theoretical modular data, I.} J. Alg. 323 (2010), 1321--1348.

\bibx{Dav2} {Davydov. A.} {Bogomolov multiplier, double class-preserving automorphisms, and modular invariants for orbifolds.} {J. Math. Phys.} {55} (2014), 092305.

\bibx{Dav3} {Davydov. A., Simmons, D.} {On Lagrangian algebras in group-theoretical braided fusion categories.} J. Alg. 471 (2017), 149--175. 


\bibx{DPR} {Dijkgraaf, R. Pasquier, V., Roche, P.} {Quasi-quantum groups related to orbifold models.}Nucl.
Phys. B. Proc. Suppl. 18B (1990), 60--72.

\bibx{DW} {Dijkgraaf, R., Witten, E.}  {Topological gauge theories and group cohomology.} {Commun. Math.
Phys.} 129 (1990), 393--429.


\bibx{DoMa} {Dong, C., Mason, G.} {Quantum Galois theory for compact Lie groups.}
{J. Algebra} {214} (1999), 92--102. 

\bibx{DRX} {Dong, C.,  Ren,  L., Xu, F.} {On orbifold theory.}  {Adv. Math.}
{321} (2017), 1--30.


\bibx{EMS} {van Ekeren, J., M\"oller, S., Scheithauer, N. R.} {Construction and classification of holomorphic vertex operator algebras.}  {J. Reine Angew. Math.} (to appear); arXiv:1507.08142.



\bibx{EGNO} {Etingof, P., Gelaki, S., Nikshych, D., Ostrik, V.} {Tensor categories.} Providence: American Math. Soc., 2015.



\bibx{Ev} {Evans, D. E.} {Twisted K-theory and modular invariants: I Quantum
doubles of finite groups.} In: Bratteli, O., Neshveyev, S., Skau, C. (eds.)
Operator Algebras: The Abel Symposium 2004. Berlin-Heidelberg: Springer, 2006,
pp. 117--144.

\bibx{EG1} {Evans, D. E., Gannon, T.}
{Modular invariants and twisted equivariant K-theory.} {Commun. Number Th. 
Phys.} {3} (2009), 209--296.


\bibx{EG2} {Evans, D. E., Gannon, T.}
{The exoticness and realisability of twisted Haagerup-Izumi modular data.} 
Commun. Math. Phys. {307} (2011), 463--512. 

\bibx{EG3} {Evans, D. E., Gannon, T.}
{Modular invariants and twisted equivariant K-theory II: Dynkin diagram symmetries.} {J. K-Theory} 12 (2013), 273--330.


\bibx{EG7} {Evans, D. E., Gannon, T.}
{Modular invariants and KK-theory I: loop groups.} (in preparation).


\bibx{EG8} {Evans, D. E., Gannon, T.}
{Modular invariants and KK-theory II: finite groups the full story.} (in preparation).




\bibx{EP1}
{Evans, D. E., Pinto, P. T.}
{Subfactor realisation of modular invariants.} {\CMP} {237} 
{(2003)}, {309--363}.



\bibx{FFSS} {Fjelstad, J., Fuchs, J., Schweigert, C., Stigner, C.} {Partition functions, mapping class groups and Drinfeld doubles.} Symmetries and Groups in Contemporary Physics,
(Tianjin, 2012), World Scientific, pp. 405--410.


\bibx{FHT} {Freed, D. S., Hopkins, M. J., Teleman, C.} 
{Twisted equivariant K-theory
    with complex coefficients.} {J. Topol. 1 (2008), 16--44.}

\bibx{FHTi} {Freed, D. S., Hopkins, M. J., Teleman, C.}  {Loop groups and twisted K-theory I.} {  J. Topol.} 4 (2011),  737--798.

\bibx{FFRS} {Fro\"hlich, J., Fuchs, J., Runkel, I.,  Schweigert, C.} {Correspondences of ribbon categories.}  {Adv. Math.} 199 (2006), 192--329.


\bibx{GK} {Gem\"unden, T., Keller, C. A.} {Orbifolds of lattice vertex operator algebras at $d=48$ and $d=72$.} arXiv:1802.10581.



\bibx{GMN} {Goff, C., Mason, G., Ng,  S.-H.} {On the gauge equivalence of twisted quantum doubles of elementary
abelian and extra-special 2-groups.} {J. Algebra} 312 (2007),  849--875.


\bibx{HJPR} {Habegger, N., Jones, V., Pino Ortiz, O., Ratcliffe, J.} {Relative cohomology of groups.} {Comment. Math. Helvetici} {59} (1984), 149--164.

\bibx{HMT} {Hanaki, A., Miyamoto, M., Tambara, D.} {Quantum Galois theory for finite groups.} {Duke Math. J.} {97} (1999), 541--544.


\bibx{HWY} {Huang, H.-L., Wan, Z., Ye, Y.} {Explicit cocycle formulas on finite abelian groups with applications to braided linear Gr-categories and Dijkgraaf-Witten invariants;} arXiv:1703.03266. 

\bibitem{Hua}  {Huang, Y.-Z.} {Rigidity and modularity of vertex tensor categories.} {Commun. Contemp. Math.} 10 (2008), 871--911. 


\bibx{HKL} {Huang, Y.-Z., Kirillov, A. Jr, Lepowsky, J.} {Braided tensor categories and
extensions of vertex operator algebras.} Commun. Math. Phys.
337 (2015), 1143--1159. 

\bibx{Hugh} {Hughes, N. J. S.} {The use of bilinear mappings in the classification of groups of class 2.} {Proc. AMS} 2 (1951), 742--747.


\bibx{JF} {Johnson-Freyd, T.} {The moonshine anomaly.} {Commun.\ Math.\ Phys.} (to appear); arXiv: 1707.08388v2.

\bibx{Jon1} {Jones, V. F. R.}  {An invariant for group actions.} {Alg\`ebres d'op\'erateurs (S\'em., Les Plans-sur-Bex, 1978),} pp. 237--253, {Lecture Notes in Math., 725, Springer, Berlin, 1979.} 

\bibx{Jon2} {Jones, V. F. R.}  {Actions of finite groups on the hyperfinite type II$_1$ factor.} {Mem. Amer. Math. Soc.} 28 (1980), no. 237.


\bibx{KLX} {Kac, V. G., Longo, R., Xu, F.} {Solitons in affine and permutation orbifolds.} {Commun. Math. Phys.} {253 (2005), 723--764.}

\bibx{Karp} {Karpilovsky, G.}  {Group Representations. Vol. 2.} North-Holland Mathematics Studies, 177. Amsterdam: North-Holland Publishing Co., 1993.


\bibx{Kaw} {Kawahigashi, Y.} {Conformal field theory, tensor categories and operator algebras.} {J.\ Phys. A 48} (2015) no.30, 303001. 

\bibx{KLM} {Kawahigashi, Y., Longo, R., M\"uger, M.} {Multi-interval subfactors and modularity of representations in conformal field theory.} {Commun. Math. Phys.} 219 (2001),  631--669.



\bibx{Kir} {Kirillov, A. Jr.}
{Modular categories and orbifold models.}
Commun. Math. Phys. 229 (2002), 309--335. 

\bibx{KO} {Kirillov, A. Jr., Ostrik, V.} {On a $q$-analogue of the McKay correspondence and the ADE classification of sl(2) conformal field theories.} {Adv. Math.} 171 (2002), 183--227.


\bibx{KMY} {Kosaki, H., Munemasa, A., Yamagami, S.} {On fusion algebras associated to finite group actions.}
Pac. J. Math. 177 (1977), 269--290.


\bibx{LL} {Lepowsky, J.,  Li, H.} {Introduction to
Vertex Operator Algebras and Their Representations.} (Birkh\"auser, Boston
2004).

\bibx{LR}
{Longo, R.,  Rehren, J.-H.} {Nets of subfactors.} {Rev. Math. Phys.} {7} (1995), 567--597.

\bibx{LX} {Longo, R., Xu, F.} {Topological sectors and a dichotomy in conformal field theory.} {Commun. Math. Phys.} 251 (2004) 321--364.

\bibx{Maj} {Majid, S.} {Quantum double for quasi-Hopf algebras.} Lett. Math. Phys. 45 (1998), 1--9.

\bibx{MN} {Mason, G., Ng, S.-H.} {Generalized twisted quantum doubles of a finite group and rational orbifolds.} {arXiv:1703.06489.}

\bibx{Miy} {Miyamoto, M.} {$C_2$-cofiniteness of cyclic-orbifold models.} {Commun. Math. Phys.} 335(3) (2015), 1279--1286.

\bibx{Mug} {M\"uger, M.} {Conformal field theory and Doplicher-Roberts reconstruction.} {In: Mathematical Physics
in Mathematics and Physics: Quantum and Operator Algebraic Aspects,} R. Longo (ed.), Fields Inst.
Commun. 20, pp.297--319 (2001).

\bibx{Mug1} {M\"uger, M.} {On the structure of braided crossed $G$-categories.} Appendix to V.Turaev, Homotopy quantum field theory, volume 10 of EMS Tracts in Mathematics. European Mathematical Society, Z\"urich, 2010.


\bibx{NN} {Naidu, D., Nikshych, D.} {Lagrangian subcategories and braided tensor equivalences of 
twisted quantum doubles of finite groups.} {Commun. Math. Phys.} 279 (2008) 845--872.


\bibx{Nak} {Nakaoka, M.} {Decomposition theorem for homology groups of symmetric groups.} {Ann. Math.} 71 (1960), 16--42.


\bibx{O}
{Ocneanu, A.} {Paths on Coxeter diagrams: From Platonic solids and
singularities to minimal models and subfactors.} (Notes recorded by S.\ Goto).
In: Rajarama Bhat, B.V.\ et al.\ (eds.) Lectures on Operator Theory.
Providence: American Mathematical Society, 2000, pp. 243--323. 



\bibx{Op} {Opolka, H.} {Group extensions and cohomology groups II.} {J. Alg.} {169} (1994), 328--331. 

\bibx{Ost1} {Ostrik, V.} {Module categories, weak Hopf algebras and modular invariants.} Transform. Groups
8 (2003), 177--206.

\bibx{Ost} {Ostrik, V.} {Module categories for quantum doubles of finite groups.}
Int. Math. Res. Notices 27 (2003), 1507--1520.


\bibx{RSW} {Raeburn, I., Sims, A., Williams, D. P.} {Twisted actions and obstructions in group cohomology.} {In: $C^*$-algebras (M\"unster, 1999),} Cuntz J., Echterhoff S. (eds), pp.161--181, Springer, Berlin, 2000. 



\bibx{xu}{Xu, F.}
{New braided endomorphisms from conformal inclusions.} \CMP{} {192},
{349--403} {(1998)}.



\bibx{Xu} {Xu, F.}
{Algebraic orbifold conformal field theories.}  
Proc. Natl. Acad. Sci. USA 97, no. 26  (2000), 14069--14073. 


\end{scriptsize}

\end{thebibliography}
\end{document}